\documentclass[11pt,reqno]{amsart}

\usepackage[a4paper,margin=1.1in]{geometry}
\usepackage{amsmath,amssymb,amsthm,mathtools,mathrsfs}
\usepackage{enumitem}
\usepackage{microtype}
\usepackage[hypertexnames=false,hidelinks]{hyperref}
\usepackage{cleveref}

\numberwithin{equation}{section}
\newtheorem{theorem}{Theorem}[section]
\newtheorem{proposition}[theorem]{Proposition}
\newtheorem{lemma}[theorem]{Lemma}
\newtheorem{corollary}[theorem]{Corollary}
\theoremstyle{remark}
\newtheorem{remark}[theorem]{Remark}

\newcommand{\tr}{\operatorname{tr}}
\newcommand{\dist}{\operatorname{dist}}
\newcommand{\II}{\mathrm{II}}
\newcommand{\doi}[1]{\href{https://doi.org/#1}{\nolinkurl{doi:#1}}}
\hypersetup{pdftitle={Power Concavity and Log-Concavity for Large-Parameter Robin 2-Hessian Problems},pdfauthor={Shujun Shi},pdfsubject={Robin stability and concavity for the 2-Hessian operator},pdfkeywords={2-Hessian, Robin boundary condition, power concavity, log-concavity, eigenfunction}}

\title[Large-parameter Robin 2-Hessian concavity]{Power Concavity and Log-Concavity for Large-Parameter Robin \(2\)-Hessian Problems}
\author{Shujun Shi}
\address{School of Mathematical Sciences, Harbin Normal University,
	Harbin 150025, Heilongjiang Province, China}
\email{shjshi@hrbnu.edu.cn}
\author{Shenghui Sun}
\address{School of Mathematical Sciences, Harbin Normal University,
	Harbin 150025, Heilongjiang Province, China}
\email{17388005715@163.com}
\author{Yuzhou Zhang}
\address{School of Mathematical Sciences, Harbin Normal University,
	Harbin 150025, Heilongjiang Province, China}
\email{yuzhouzhang@stu.hrbnu.edu.cn}
\date{}

\begin{document}
\raggedbottom

\begin{abstract}
Let $\Omega\subset\mathbb R^n$, $n\ge3$, be a bounded uniformly convex
$C^{3,\gamma}$ domain, where $0<\gamma<1$.  We study the limit as the
positive Robin parameter $\alpha\to\infty$.  For the $2$-Hessian torsion problem, let $U$ be
the admissible Dirichlet solution.  For all sufficiently large $\alpha$,
we obtain a unique admissible Robin solution $u_\alpha$ and establish
the optimal-order estimate
\[
 \|u_\alpha-U\|_{C^{2,\gamma}(\overline\Omega)}
 \le C\alpha^{-1}.
\]
Moreover, $-\sqrt{-u_\alpha}$ has a uniform positive lower bound on
its Hessian, and hence $\sqrt{-u_\alpha}$ is strictly concave.  For the normalized eigenvalue problem, let
$(\lambda_D,\varphi_D)$ be the normalized negative admissible Dirichlet
eigenpair, and assume $\varphi_D\in C^{2,\gamma}(\overline\Omega)$.  For all sufficiently large $\alpha$, there is a unique normalized
negative admissible Robin eigenpair $(\lambda_\alpha,\varphi_\alpha)$,
which satisfies
\[
 \|\varphi_\alpha-\varphi_D\|_{C^{2,\gamma}(\overline\Omega)}
   +|\lambda_\alpha-\lambda_D|\le C\alpha^{-1}.
\]
The Hessian of $-\log(-\varphi_\alpha)$ also has a uniform positive
lower bound.
The main analytic tool is a uniform H\"older estimate for linear
Robin problems, with an $\alpha^{-1}$ gain for $C^{2,\gamma}$ boundary data.  A normalized augmented inverse removes the scaling
kernel in the eigenvalue problem.  A Schur-complement argument combines the $C^2$ convergence with
boundary geometry to prove convexity up to the boundary and anisotropic
Hessian lower bounds at the scale
$\dist(x,\partial\Omega)+\alpha^{-1}$.  The convergence rate $O(\alpha^{-1})$ for the Robin solution
functions is optimal even in $C^0$, as follows from their boundary
conditions.
\end{abstract}
\maketitle

\section{Introduction}

Convexity and concavity of solutions and their natural nonlinear
transforms have long been studied in elliptic PDEs.  For the Laplacian,
the square root of the torsion function is concave on convex domains,
and the first Dirichlet eigenfunction is log-concave; see
\cite{BL,Ken,ML}.  Robin boundary conditions alter the behavior of these transforms near
the boundary.  As the positive Robin parameter tends to infinity, the
solutions approach their Dirichlet counterparts.  However, the Hessians
of the square-root and logarithmic transforms contain negative powers
of the solutions, which become singular near the Dirichlet boundary.
Thus interior $C^2$ convergence alone does not directly yield
concavity estimates uniform up to the boundary.

For the Laplacian, large-parameter Robin concavity has been studied in
\cite{CF,YZ}.  In particular, Ye--Zhang proved quantitative $C^2$
convergence on $C^{3,1}$ uniformly convex domains and used a boundary
argument to derive the corresponding concavity.  Robin concavity need not hold for every parameter value, as shown
by counterexamples for the linear ground state in \cite{ACH}.

For the $2$-Hessian operator, the corresponding Dirichlet convexity
results have been established.  In dimension three, results for the torsion transform
$-\sqrt{-U}$ were obtained by Ma--Xu \cite[Theorem~2]{MX} and Salani
\cite[Theorem~1.2]{Salani}.  For the logarithmic eigenfunction transform
$-\log(-\varphi_D)$, see Liu--Ma--Xu
\cite[Theorem~1.1]{LMX} and Salani \cite[Theorem~1.1]{Salani}.  The local
constant-rank results used to pass from a positive semidefinite Hessian
to strict convexity appear in \cite[Lemma~1]{MX} and
\cite[Theorem~2.1]{LMX}.  For
$n\ge4$, the corresponding Dirichlet convexity and constant-rank results
are available in Li--Ma--Qiu--Salani \cite{LMQS}.

For fixed Robin parameters, solvability of Hessian equations with
oblique boundary conditions is covered by the theory of
Lions--Trudinger--Urbas \cite{LTU} and Ma--Qiu \cite{MQ}.  For each fixed positive Robin parameter, a spatial dilation
reduces the torsion problem to the coefficient-one boundary condition
treated by Ma--Qiu on a rescaled uniformly convex domain.
These fixed-parameter results do not yield estimates uniform as
$\alpha\to\infty$.  Nor do they provide a spectral inverse that
incorporates the normalization.

We study these large-parameter Robin problems for the nonlinear
$2$-Hessian operator.  Compared with the linear theory,
three issues have to be addressed.  First, the linearization depends on the Dirichlet solution, and its
Robin estimates must be uniform as $\alpha\to\infty$.  Second, the Dirichlet linearization of the eigenvalue equation has
the eigenfunction itself as a scaling mode in its kernel.  The inverse
must therefore account for both the eigenvalue variation and the
normalization.  Third, the singular
square-root and logarithmic Hessians cannot be controlled near
$\partial\Omega$ by interior $C^2$ convergence alone.  We address these difficulties through a uniform Robin resolvent, an
augmented spectral inverse, and a boundary Schur-complement estimate,
respectively.

For a symmetric matrix $A$, write
\[
\sigma_1(A)=\operatorname{tr}A,\qquad
\sigma_2(A)=\frac12\bigl((\operatorname{tr}A)^2-\operatorname{tr}(A^2)\bigr),
\qquad
T_1(A)=\sigma_1(A)I-A,
\]
and let
\[
\Gamma_2
=\{A:\sigma_1(A)>0,\ \sigma_2(A)>0\}.
\]
Let $\nu$ denote the outward unit normal.  The Robin torsion problem is
\begin{equation}\label{eq:RT}
 \begin{cases}
 \sigma_2(D^2u_\alpha)=1&\text{in }\Omega,\\
 (u_\alpha)_\nu+\alpha u_\alpha=0&\text{on }\partial\Omega,\\
 u_\alpha<0,\qquad D^2u_\alpha\in\Gamma_2,
 \end{cases}
 \tag{$R_\alpha$}
\end{equation}
and its Dirichlet limit is
\begin{equation}\label{eq:DT}
 \sigma_2(D^2U)=1\quad\text{in }\Omega,\qquad U=0\quad\text{on }\partial\Omega.
 \tag{$D$}
\end{equation}
Here $U$ is the negative admissible solution of the Dirichlet problem.
The normalized Robin eigenvalue problem is
\begin{equation}\label{eq:RE}
 \begin{cases}
 \sigma_2(D^2\varphi_\alpha)=\lambda_\alpha(-\varphi_\alpha)^2&\text{in }\Omega,\\
 (\varphi_\alpha)_\nu+\alpha\varphi_\alpha=0&\text{on }\partial\Omega,\\
 \displaystyle\int_\Omega(-\varphi_\alpha)^3=1,\quad
 \varphi_\alpha<0,\quad D^2\varphi_\alpha\in\Gamma_2.
 \end{cases}
 \tag{$E_\alpha$}
\end{equation}
The corresponding Dirichlet eigenpair $(\lambda_D,\varphi_D)$ satisfies
\[
 \sigma_2(D^2\varphi_D)=\lambda_D(-\varphi_D)^2
 \quad\text{in }\Omega,\qquad
 \varphi_D=0\quad\text{on }\partial\Omega,
\]
where $\lambda_D>0$, $\varphi_D<0$ and
$D^2\varphi_D\in\Gamma_2$ in $\Omega$, with
$\int_\Omega(-\varphi_D)^3\,dx=1$.

Throughout this paper, a classical Robin solution, including the
eigenfunction in a Robin eigenpair, belongs to $C^2(\overline\Omega)$.
It satisfies the equation in $\Omega$ and the Robin boundary condition
pointwise on $\partial\Omega$.  Admissibility means
$D^2u(x)\in\Gamma_2$ for every $x\in\Omega$ (with $u$ replaced by
$\varphi$ for an eigenfunction).

The spectral result requires an additional Dirichlet boundary
regularity assumption:
\begin{equation}\label{eq:Eg}
 \tag{$E_\gamma$}
 \text{the normalized negative admissible Dirichlet eigenpair exists and }
 \varphi_D\in C^{2,\gamma}(\overline\Omega).
\end{equation}
On smooth domains, existence and global $C^{1,1}$ regularity follow
from the classical $k$-Hessian eigenvalue theory; see Wang
\cite{Wang,Wang09}.  For domains of finite regularity, the existence of the Dirichlet
eigenpair is explicitly included in \eqref{eq:Eg}.
Lemma~\ref{lem:dirichlet-spectral} derives uniform ellipticity, an
additional derivative, and the augmented Dirichlet inverse directly
from this assumption, without using logarithmic convexity.  The condition holds on every ball by
Proposition~\ref{prop:Eg-ball}; further comments are given in
Remark~\ref{rem:Eg-domains}.

We can now state the main results.

\begin{theorem}[Large-parameter torsion solutions and power concavity]\label{thm:torsion-main}
Let $n\ge3$, $0<\gamma<1$, and let $\Omega\subset\mathbb R^n$ be a
bounded uniformly convex $C^{3,\gamma}$ domain.  There exist
$\alpha_T,C_T,c_T>0$ such that, for every $\alpha\ge\alpha_T$, problem
\eqref{eq:RT} has a unique admissible classical solution and
\begin{equation}\label{eq:torsion-rate-main}
 \|u_\alpha-U\|_{C^{2,\gamma}(\overline\Omega)}\le C_T\alpha^{-1}.
\end{equation}
Moreover,
\begin{equation}\label{eq:torsion-concavity-main}
 D^2[-\sqrt{-u_\alpha}]\ge c_TI\qquad\text{in }\Omega.
\end{equation}
In particular, $\sqrt{-u_\alpha}$ is strictly concave for every
$\alpha\ge\alpha_T$.
\end{theorem}

\begin{theorem}[Large-parameter eigenpairs and log-concavity]\label{thm:eigen-main}
Let $n\ge3$, $0<\gamma<1$, and let $\Omega$ be a
bounded uniformly convex $C^{3,\gamma}$ domain.  Assume \eqref{eq:Eg}.  There exist
$\alpha_E,C_E,c_E>0$ such that, for every $\alpha\ge\alpha_E$, problem
\eqref{eq:RE} has a unique normalized negative admissible classical
eigenpair and
\begin{equation}\label{eq:eigen-rate-main}
 \|\varphi_\alpha-\varphi_D\|_{C^{2,\gamma}(\overline\Omega)}
  +|\lambda_\alpha-\lambda_D|\le C_E\alpha^{-1}.
\end{equation}
Moreover,
\begin{equation}\label{eq:eigen-concavity-main}
 D^2[-\log(-\varphi_\alpha)]\ge c_EI\qquad\text{in }\Omega.
\end{equation}
Hence, $-\varphi_\alpha$ is strictly log-concave for every
$\alpha\ge\alpha_E$.
\end{theorem}

The proof also identifies the boundary scale and provides anisotropic
lower bounds for the transformed Hessians.

\begin{corollary}[Boundary scale and anisotropic convexity]\label{cor:boundary-scale}
After increasing $\alpha_T$, and also $\alpha_E$ under~\eqref{eq:Eg},
if necessary, there are constants $\delta,c,C>0$ for which the following
estimates hold for every $\alpha\ge\alpha_T$ in the torsion case and,
under \eqref{eq:Eg}, for every $\alpha\ge\alpha_E$ in the spectral case.
The constants $c,C$ may differ between the two problems.  We take
$\delta$ smaller than the tubular radius of $\partial\Omega$.
For $x\in\overline\Omega$ with $d:=\dist(x,\partial\Omega)<\delta$,
there is a unique point $y\in\partial\Omega$ such that
$x=y-d\nu(y)$.  Then
\begin{equation}\label{eq:intro-boundary-size-torsion}
 c(d+\alpha^{-1})\le -u_\alpha(x)
 \le C(d+\alpha^{-1}),
\end{equation}
and, for every $\zeta\in T_y\partial\Omega$ and $s\in\mathbb R$,
\begin{equation}\label{eq:intro-anisotropic-torsion}
 D^2[-\sqrt{-u_\alpha}][\zeta+s\nu,\zeta+s\nu]
 \ge c\Bigl((d+\alpha^{-1})^{-1/2}|\zeta|^2
 +(d+\alpha^{-1})^{-3/2}s^2\Bigr).
\end{equation}
Under \eqref{eq:Eg}, the corresponding eigenfunction satisfies
\begin{equation}\label{eq:intro-boundary-size-eigen}
 c(d+\alpha^{-1})\le -\varphi_\alpha(x)
 \le C(d+\alpha^{-1}),
\end{equation}
and
\begin{equation}\label{eq:intro-anisotropic-eigen}
 D^2[-\log(-\varphi_\alpha)][\zeta+s\nu,\zeta+s\nu]
 \ge c\Bigl((d+\alpha^{-1})^{-1}|\zeta|^2
 +(d+\alpha^{-1})^{-2}s^2\Bigr).
\end{equation}
\end{corollary}

For the eigenvalue problem, we first establish uniform ellipticity
and invertibility of the augmented Dirichlet linearization under
\eqref{eq:Eg}, without using log-concavity.  The implicit function
theorem then provides smooth-domain eigenpairs converging to the
prescribed Dirichlet eigenpair.  This allows the known convexity
theorems to pass to $\Omega$.  Uniform Robin estimates yield the
$\alpha^{-1}$ bounds, and a boundary Schur-complement argument proves
positivity of the transformed Robin Hessians.  Corollary~\ref{cor:boundary-scale}
identifies $\alpha^{-1}$ as both the convergence rate and the boundary
scale of these Hessians.

The paper is organized as follows. Section~\ref{sec:dirichlet-input} establishes the required Dirichlet
properties and linear estimates.  Section~\ref{sec:branches} constructs
the large-$\alpha$ solutions, and Section~\ref{sec:concavity} proves
the transformed-Hessian bounds and the main theorems.  The final section
concerns geometric consequences, the spectral boundary hypothesis on
balls, and the optimality and uniformity of the estimates.

\section{Dirichlet estimates and linear Robin theory}\label{sec:dirichlet-input}

\subsection{Algebraic and boundary identities}
We use the following notation throughout.  Let
$\operatorname{Sym}(n)$ denote the space of real symmetric
$n\times n$ matrices.  For $A,B\in\operatorname{Sym}(n)$, set
\[
A:B=\operatorname{tr}(AB),
\qquad
|A|=(A:A)^{1/2}.
\]
Here $:$ denotes the Frobenius inner product; $|\cdot|$ denotes the
Frobenius norm on matrices and the Euclidean norm on vectors.  We write
\[
\sigma_1(A)=\operatorname{tr}A,
\qquad
\sigma_2(A)
=\frac12\bigl((\operatorname{tr}A)^2-\operatorname{tr}(A^2)\bigr),
\]
and
\[
T_1(A)=\sigma_1(A)I-A,
\qquad
\Gamma_2
=\bigl\{A\in\operatorname{Sym}(n):
\sigma_1(A)>0,\ \sigma_2(A)>0\bigr\}.
\]
Here $I$ is the identity matrix.  The inequalities
$A\ge0$ and $A>0$ mean that $A$ is positive semidefinite and positive
definite, respectively, and $A\ge cI$ means $A-cI\ge0$.  For
$\xi,\eta\in\mathbb R^n$, we use
\[
 A[\xi,\eta]=\xi^TA\eta,
 \qquad
 A[\xi]=A[\xi,\xi],
 \qquad
 (\xi\otimes\eta)_{ij}=\xi_i\eta_j.
\]
The notation $D^2v>0$ means that the Hessian is positive definite
at every point, which implies strict convexity on a convex domain.
The stronger inequality $D^2v\ge cI$ with $c>0$ gives a uniform
Hessian lower bound.
Repeated indices are summed, with $i,j$ ranging over
$1,\dots,n$ and tangential indices $a,b$ over $1,\dots,n-1$.  We write $[h]_{C^\gamma}$ for the H\"older seminorm.
We denote by $\nu$ the outward unit normal to $\partial\Omega$, set
\[
 f_\nu=Df\cdot\nu,
 \qquad
 d(x)=\operatorname{dist}(x,\partial\Omega),
\]
and identify $T_x\partial\Omega$ with $\nu(x)^\perp$.

For $A,B\in\operatorname{Sym}(n)$,
\begin{align}
 D\sigma_2(A)[B]&=T_1(A):B,\label{eq:D-sigma}\\
 \sigma_2(A+B)&=\sigma_2(A)+T_1(A):B+\sigma_2(B),\label{eq:quad}\\
 \sigma_2(A)-\sigma_2(B)&=
 T_1\!\left(\frac{A+B}{2}\right):(A-B).\label{eq:difference}
\end{align}
If $A\in\Gamma_2$, then $T_1(A)>0$; see \cite{Garding}.  At a
boundary point $x\in\partial\Omega$, relative to the decomposition
\[
 \mathbb R^n=T_x\partial\Omega\oplus\mathbb R\nu(x),
\]
write
\[
 A=\begin{pmatrix}A_T&b\\b^T&r\end{pmatrix},
\]
where $A_T$ is the tangential block, $b\in T_x\partial\Omega$, and
$r=A[\nu,\nu]$.  Then
\begin{equation}\label{eq:block}
 \sigma_2(A)=\sigma_2(A_T)+r\,\tr A_T-|b|^2,
 \qquad T_1(A)[\nu,\nu]=\tr A_T.
\end{equation}
At $x\in\partial\Omega$, let
$\II_x(\tau,\zeta)=\langle D_\tau\nu,\zeta\rangle$ be the second
fundamental form.  With this sign convention, uniform convexity means
\[
 \II_x\ge\kappa_0 I\quad\text{on }T_x\partial\Omega
\]
for some $\kappa_0>0$ independent of $x$.  If $f=0$ on
$\partial\Omega$, then, at each $x\in\partial\Omega$,
\begin{equation}\label{eq:boundary-hessian}
 D^2f(\tau,\zeta)=f_\nu\II_x(\tau,\zeta)
 \qquad(\tau,\zeta\in T_x\partial\Omega).
\end{equation}

\subsection{Properties of the Dirichlet solutions}
We first establish the boundary ellipticity and augmented Dirichlet
invertibility needed for a strong approximation of the spectral pair.
These arguments use \eqref{eq:Eg}, but not convexity of its logarithmic
transform.  The known smooth-domain convexity results will be applied
only after that approximation has been proved.  The following algebraic
and full-rank lemmas will also be used for the Dirichlet torsion problem.

\begin{lemma}[Boundary Newton tensor]\label{lem:boundary-newton}
Let $n\ge3$ and
\[
 B=\begin{pmatrix}A&b\\ b^T&r\end{pmatrix},
 \qquad A\in\operatorname{Sym}(n-1),\qquad A>0,
\]
and assume $\sigma_2(B)=0$.  Then $T_1(B)>0$.  More precisely, if
$0<a_*\le a^*$, $a_*I\le A\le a^*I$, and $|b|\le M$, the eigenvalues of $T_1(B)$ are bounded above and bounded away from
zero by constants depending only on $n,a_*,a^*,M$.
\end{lemma}
\begin{proof}
Set $H=\tr A$.  Formula \eqref{eq:block} yields
\[
 r=\frac{|b|^2-\sigma_2(A)}H.
\]
The lower-right entry of $T_1(B)$ equals $H$, and the corresponding
Schur complement is
\begin{equation}\label{eq:newton-schur}
 S=HI-A-\frac{\sigma_2(A)}H I+\frac{|b|^2I-b\otimes b}{H}.
\end{equation}
The last term is positive semidefinite.  Let $a_i>0$ be the eigenvalues of $A$, and set
$S_i=\sum_{j\ne i}a_j$.  Then
\[
 H(H-a_i)-\sigma_2(A)
 =S_i^2-\sum_{\substack{j<k\\j,k\ne i}}a_ja_k
 =\sum_{j\ne i}a_j^2+
   \sum_{\substack{j<k\\j,k\ne i}}a_ja_k>0.
\]
It follows that
$HI-A-\frac{\sigma_2(A)}H I$ in \eqref{eq:newton-schur} is positive
definite.  Its eigenvalues are bounded below by $a_*^2/H$.  Since
$(n-1)a_*\le H\le(n-1)a^*$ and $|b|\le M$, the Schur-complement
criterion also supplies uniform upper and lower eigenvalue bounds
for $T_1(B)$ with the stated dependence.
\end{proof}

\begin{lemma}[Full rank from the gradient image]\label{lem:full-rank}
Let $D\subset\mathbb R^n$ be a bounded domain, and let
$v\in C^2(D)$ satisfy $D^2v\ge0$ and have constant Hessian rank.  If
either
\begin{enumerate}[label=(\alph*),leftmargin=*]
\item $v\in C(\overline D)$, $v<0$ in $D$, and $v=0$ on $\partial D$; or
\item $v(x)\to+\infty$ as $x\to\partial D$,
\end{enumerate}
then $D^2v>0$ in $D$.
\end{lemma}
\begin{proof}
In case (a), choose $x_0\in D$ and set $a=-v(x_0)>0$.  For every
$|p|<a/(2\operatorname{diam}D)$, the function $v-p\cdot x$ attains an
interior minimum at $x$, where $Dv(x)=p$.  Thus $Dv(D)$ contains an open ball.  In case
(b), $v-p\cdot x$ has an interior minimum for every $p$, so
$Dv(D)=\mathbb R^n$.  If the constant Hessian rank were $r<n$, the constant-rank theorem
for the $C^1$ map $Dv$ would cover its image by countably many
$r$-dimensional submanifolds.  Such a set cannot contain an open ball,
contrary to both cases.
\end{proof}

\begin{lemma}[Dirichlet spectral regularity and augmented invertibility]
\label{lem:dirichlet-spectral}
Let $n\ge3$, $0<\gamma<1$, and let $\Omega$ be a bounded uniformly
convex $C^{3,\gamma}$ domain.  Assume \eqref{eq:Eg}.  Then
$\varphi_D\in C^{3,\gamma}(\overline\Omega)\cap C^\infty(\Omega)$,
$q_E:=(\varphi_D)_\nu>0$ on $\partial\Omega$, and there are
$\vartheta_E,\Theta_E>0$ such that
\begin{equation}\label{eq:spectral-base-ellipticity}
 \vartheta_E I\le a_E:=T_1(D^2\varphi_D)\le\Theta_E I
 \quad\text{on }\overline\Omega.
\end{equation}
Define
\begin{equation}\label{eq:linearizations}
 L_Ev=a_E:D^2v,\qquad
 M_Ev=L_Ev-2\lambda_D\varphi_Dv,\qquad
 \mathfrak n(v)=\int_\Omega\varphi_D^2v\,dx.
\end{equation}
For every $0<\eta\le\gamma$, set
\[
 X_{0,\eta}=\{v\in C^{2,\eta}(\overline\Omega):
                      v|_{\partial\Omega}=0\}.
\]
The augmented Dirichlet operator
\begin{equation}\label{eq:dir-aug-map}
 \mathcal B_D(v,\ell)
   =\bigl(M_Ev-\ell\varphi_D^2,\mathfrak n(v)\bigr)
\end{equation}
is a bounded isomorphism from $X_{0,\eta}\times\mathbb R$ onto
$C^{0,\eta}(\overline\Omega)\times\mathbb R$.  In particular,
\begin{equation}\label{eq:dir-aug-est}
 \|v\|_{C^{2,\eta}}+|\ell|
 \le C_\eta\bigl(\|M_Ev-\ell\varphi_D^2\|_{C^{0,\eta}}
                       +|\mathfrak n(v)|\bigr).
\end{equation}
None of these conclusions requires convexity of
$-\log(-\varphi_D)$.
\end{lemma}
\begin{proof}
We first establish boundary ellipticity, then improve the regularity,
and finally determine the kernel and range of the Dirichlet operator.

\smallskip
\noindent\emph{Boundary ellipticity and regularity.}
By \eqref{eq:Eg}, $D^2\varphi_D$ is continuous on $\overline\Omega$.
Admissibility implies $\Delta\varphi_D>0$ in $\Omega$, so the Hopf
boundary lemma gives
\[
 q_E=(\varphi_D)_\nu>0,\qquad
 q_{E,*}:=\min_{\partial\Omega}q_E>0.
\]
At $x\in\partial\Omega$, tangential differentiation of the zero
Dirichlet trace yields
\[
 D^2\varphi_D|_{T_x\partial\Omega}=q_E\II_x,\qquad
 D^2\varphi_D(\tau,\nu)=\partial_\tau q_E.
\]
The equation extends continuously to the boundary, where
$\sigma_2(D^2\varphi_D)=0$.  Apply Lemma~\ref{lem:boundary-newton}
with $A=q_E\II_x$ and $b=\nabla_{\partial\Omega}q_E$.
It gives a positive lower bound for $T_1(D^2\varphi_D)$ on the
boundary.  This tensor is positive definite in the interior by
admissibility.  Continuity and compactness prove
\eqref{eq:spectral-base-ellipticity}.  Taking traces also gives
\begin{equation}\label{eq:reference-trace-gap}
 \Delta\varphi_D=\frac{\operatorname{tr}a_E}{n-1}
       \ge\frac{n\vartheta_E}{n-1}>0
       \quad\text{on }\overline\Omega.
\end{equation}

In a $C^{3,\gamma}$ flattened boundary chart, write the equation
as $\mathcal G(y,u,Du,D^2u)=0$ for the pullback $u$ of $\varphi_D$.
For a tangential increment $h\tau$, the difference quotient
$z_h=(u(\cdot+h\tau)-u)/h$ has zero trace on the flat boundary and
satisfies a linear equation
\[
 a_h^{ij}(z_h)_{ij}+b_h^i(z_h)_i+c_hz_h=f_h.
\]
The principal coefficients are averages of the derivatives of
$\mathcal G$ with respect to its Hessian variable along the two
interpolated jets.  Uniform ellipticity of the reference linearization
and the $C^{2,\gamma}$ regularity of $u$ give uniform ellipticity and
$C^{0,\gamma}$ bounds for $a_h$.  The lower-order coefficients and
$f_h$ are controlled by the $C^{3,\gamma}$ chart and the reference
solution, uniformly for small $h$.  Boundary Schauder estimates
\cite{GT} therefore control $z_h$ in $C^{2,\gamma}$ and give all
third derivatives with a tangential index.  The derivative of the
original equation with respect to the normal-normal Hessian entry
is $a_E[\nu,\nu]>0$.  Solving for that entry and differentiating
recovers the remaining third derivative.  Thus
$\varphi_D\in C^{3,\gamma}(\overline\Omega)$ and
$a_E\in C^{1,\gamma}(\overline\Omega)$; interior differentiation and
Schauder estimates give $\varphi_D\in C^\infty(\Omega)$.
Only \eqref{eq:Eg}, boundary geometry, and uniform ellipticity have
been used.

\smallskip
\noindent\emph{The Dirichlet kernel and range.}
Write $\phi=\varphi_D$ and $\omega=-\phi>0$.  The tensor $a_E$ is
symmetric and divergence-free:
\[
 \partial_i(a_E)^{ij}
       =\partial_j\Delta\phi-\partial_i\phi_{ij}=0.
\]
Euler's identity for $\sigma_2$ gives $M_E\phi=0$.  For
$v\in X_{0,\eta}$, direct calculation yields
\begin{equation}\label{eq:ground-state-div}
 a_E\nabla v\cdot\nabla v+2\lambda_D\phi v^2
 -\omega^2a_E\nabla(v/\omega)\cdot\nabla(v/\omega)
 =\operatorname{div}\left(\frac{v^2}{\omega}a_E\nabla\omega\right).
\end{equation}
Near the boundary, $v=d\widetilde v$ and
$\omega=d\widetilde\omega$, where $\widetilde\omega$ is bounded away
from zero.  The quotient $v/\omega$ has bounded first derivatives,
and $v^2/\omega=O(d)$.  Integration over $\{d>\varepsilon\}$ and
passage to $\varepsilon\downarrow0$ therefore give
\begin{equation}\label{eq:ground-state-expanded}
 -\int_\Omega vM_Ev\,dx
 =\int_\Omega \omega^2a_E\nabla(v/\omega)\cdot\nabla(v/\omega)\,dx
 \ge0.
\end{equation}
If $M_Ev=0$, the quotient $v/\omega$ is constant on the connected
domain $\Omega$.  Thus the Dirichlet kernel is
$\operatorname{span}\{\phi\}$.

The Dirichlet operator $L_E:X_{0,\eta}\to C^{0,\eta}$ is an
isomorphism by Schauder theory and the maximum principle.  The
zeroth-order term in $M_E$ is compact between these spaces, so
$M_E$ is Fredholm of index zero.  Moreover,
\[
 M_Ev=\partial_i\bigl((a_E)^{ij}\partial_jv\bigr)
                         -2\lambda_D\phi v.
\]
The symmetric Dirichlet realization is self-adjoint in $L^2$.
In particular, Green's identity gives
$\int_\Omega\phi M_Ev=\int_\Omega vM_E\phi=0$ for
$v\in X_{0,\eta}$.  Since the range has codimension one, its
H\"older realization is exactly
\[
 \operatorname{Ran}M_E
   =\{f\in C^{0,\eta}(\overline\Omega):\int_\Omega\phi f\,dx=0\}.
\]
Consequently $M_Ev-\ell\phi^2=f$ with zero Dirichlet trace is
solvable precisely when
\[
 0=\int_\Omega\phi f\,dx+\ell\int_\Omega\phi^3\,dx.
\]
The normalization gives $\int_\Omega\phi^3=-1$, and hence
$\ell=\int_\Omega\phi f$.  Adding a multiple of $\phi$ fixes
$\mathfrak n(v)=s$ uniquely because $\mathfrak n(\phi)=-1$.
The bounded inverse theorem proves \eqref{eq:dir-aug-est} for every
$0<\eta\le\gamma$.
\end{proof}

\begin{lemma}[Smooth-domain approximation by the implicit function theorem]
\label{lem:spectral-approx}
Under the hypotheses of Lemma~\ref{lem:dirichlet-spectral}, there
exist smooth uniformly convex domains $\Omega_j$ and diffeomorphisms
$\Phi_j:\overline\Omega\to\overline{\Omega_j}$ such that
$\Phi_j\to\operatorname{id}$ in $C^{3,\eta}$ for every
$0<\eta<\gamma$.  If $(\lambda_j,\varphi_j)$ is the normalized
negative admissible Dirichlet eigenpair on $\Omega_j$, then
\begin{equation}\label{eq:spectral-pullback-convergence}
 |\lambda_j-\lambda_D|
 +\|\varphi_j\circ\Phi_j-\varphi_D\|_{C^{2,\eta}(\overline\Omega)}
 \longrightarrow0\qquad(0<\eta<\gamma).
\end{equation}
In particular,
\begin{equation}\label{eq:spectral-approx}
 \lambda_j\longrightarrow\lambda_D,\qquad
 \varphi_j\longrightarrow\varphi_D
       \quad\text{in }C^2(K)
\end{equation}
for every compact set $K\Subset\Omega$.
\end{lemma}
\begin{proof}
Fix $0<\eta<\gamma$.  We pull the equation and its normalization
back to $\Omega$, apply the augmented Dirichlet inverse from
Lemma~\ref{lem:dirichlet-spectral}, and identify the resulting
solutions only on the smooth approximating domains.

\smallskip
\noindent\emph{Smooth domains and pullback of the Hessian.}
Translate so that $0\in\Omega$ and write
$\partial\Omega=\{r_\Omega(\omega)\omega:\omega\in\mathbb S^{n-1}\}$.
The radial function $r_\Omega$ is positive and $C^{3,\gamma}$.
Choose positive smooth functions $r_j$ converging to $r_\Omega$
in $C^{3,\eta}$ for every $\eta<\gamma$, with uniformly bounded
$C^{3,\gamma}$ norms.  Their radial graphs bound smooth uniformly
convex domains $\Omega_j$ for large $j$, since uniform convexity
persists under small $C^2$ perturbations.  Choose a smooth radial
cutoff $\chi$ which vanishes near $0$ and equals $1$ on a neighborhood
of the range of $r_\Omega$.  The maps
\[
 \Phi_j(r\omega)
 =\bigl(r+\chi(r)(r_j(\omega)-r_\Omega(\omega))\bigr)\omega
\]
are diffeomorphisms onto $\Omega_j$ for large $j$ and converge to the
identity in $C^{3,\eta}(\overline\Omega)$.

Let $\Phi$ lie in a sufficiently small open neighborhood of the
identity in $C^{2,\eta}(\overline\Omega;\mathbb R^n)$; the unknown
$v$ remains in the fixed zero-Dirichlet space $X_{0,\eta}$.  Set
$P_\Phi=(D\Phi)^{-1}$ and $J_\Phi=\det D\Phi>0$.
If $u=v\circ\Phi^{-1}$ on $\Phi(\Omega)$, the chain rule gives
\begin{equation}\label{eq:pullback-hessian}
 (D^2u)\circ\Phi=\mathcal H_\Phi[v]
 :=P_\Phi^T\left(D^2v-
       \sum_{k=1}^n(P_\Phi^TDv)_kD^2\Phi^k\right)P_\Phi.
\end{equation}
Define, on $X_{0,\eta}\times\mathbb R$,
\begin{equation}\label{eq:dirichlet-ift-map}
 \mathcal F(\Phi,v,\Lambda)
 =\left(\sigma_2(\mathcal H_\Phi[v])-\Lambda v^2,
           \int_\Omega(-v)^3J_\Phi\,dx-1\right).
\end{equation}
This map takes values in $C^{0,\eta}(\overline\Omega)\times\mathbb R$
and is continuously Fr\'echet differentiable.  Indeed,
\eqref{eq:pullback-hessian} involves only differentiation, products,
and matrix inversion near the identity in H\"older spaces.  The
Jacobian in \eqref{eq:dirichlet-ift-map} is essential: it imposes
the normalization on the moving domain.  We use the polynomial
$\sigma_2$, not $\sqrt{\sigma_2}$; the vanishing of the density on
the Dirichlet boundary therefore causes no differentiability problem.

\smallskip
\noindent\emph{The augmented implicit function theorem.}
At $(\operatorname{id},\varphi_D,\lambda_D)$ the map $\mathcal F$
vanishes, and its derivative in $(v,\Lambda)$ is
\begin{equation}\label{eq:dirichlet-ift-derivative}
 D_{(v,\Lambda)}\mathcal F[h,\ell]
   =\bigl(M_Eh-\ell\varphi_D^2,-3\mathfrak n(h)\bigr).
\end{equation}
By Lemma~\ref{lem:dirichlet-spectral}, this is a bounded
isomorphism.  The Banach-space implicit function theorem gives
solutions $(v_\Phi,\Lambda_\Phi)$ near
$(\varphi_D,\lambda_D)$, with $v_\Phi=0$ on $\partial\Omega$, and
\begin{equation}\label{eq:ift-local-continuity}
 \|v_\Phi-\varphi_D\|_{C^{2,\eta}}
       +|\Lambda_\Phi-\lambda_D|\longrightarrow0
       \quad\text{as }\Phi\to\operatorname{id}\text{ in }C^{2,\eta}.
\end{equation}
The use of $\eta<\gamma$ avoids requiring smooth approximation in
the endpoint H\"older norm.  No assertion about convergence in that
endpoint norm is needed.

\smallskip
\noindent\emph{Negativity, admissibility, and identification.}
The Hopf slope $q_{E,*}>0$ and the zero trace imply
$\varphi_D\le-cd$ near $\partial\Omega$.  The $C^1$ convergence in
\eqref{eq:ift-local-continuity} preserves a positive outward
normal component of $Dv_\Phi$ on a fixed boundary neighborhood.
Integrating from the zero boundary trace gives $v_\Phi<0$ there;
on the remaining compact set the same conclusion follows from
uniform convergence to $\varphi_D<0$.  Thus $v_\Phi<0$ in $\Omega$.
Also $\Lambda_\Phi>0$, and continuity of
\eqref{eq:pullback-hessian} and \eqref{eq:spectral-base-ellipticity}
gives
\[
 T_1(\mathcal H_\Phi[v_\Phi])\ge\tfrac12\vartheta_E I
                   \quad\text{on }\overline\Omega.
\]
For $u_\Phi=v_\Phi\circ\Phi^{-1}$, this implies
$\Delta u_\Phi>0$.  The equation gives
$\sigma_2(D^2u_\Phi)=\Lambda_\Phi u_\Phi^2>0$ in $\Phi(\Omega)$.
Hence $u_\Phi$ is negative and admissible in the interior, has zero
Dirichlet trace, and satisfies
$\int_{\Phi(\Omega)}(-u_\Phi)^3=1$.

Apply this construction to $\Phi=\Phi_j$.  Each $u_{\Phi_j}$ is
$C^{2,\eta}(\overline{\Omega_j})$, hence $C^{1,1}$, and interior
Schauder bootstrapping gives $C^\infty(\Omega_j)$.  The domain
$\Omega_j$ is smooth and uniformly convex, so Wang's uniqueness
theorem \cite[Theorem~6.2(ii)]{Wang09} applies in precisely its
classical solution class.  It identifies
$(\Lambda_{\Phi_j},u_{\Phi_j})$ with
$(\lambda_j,\varphi_j)$; the integral normalization fixes the
multiplicative constant.  Thus \eqref{eq:spectral-pullback-convergence}
follows from \eqref{eq:ift-local-continuity}.  Since $\eta<\gamma$
was arbitrary and the smooth-domain eigenpairs are unique, this
holds for every such $\eta$.

For $K\Subset\Omega$, we have $K\subset\Omega_j$ for large $j$.
The chain rule \eqref{eq:pullback-hessian} and
$\Phi_j^{-1}\to\operatorname{id}$ in $C^2$ give
\eqref{eq:spectral-approx} in the original Euclidean coordinates.
In particular, the argument neither invokes uniqueness for a weak
limit on $\Omega$ nor assumes convexity of the logarithmic transform
in constructing the approximations.
\end{proof}

We use these lemmas and the standard Dirichlet theory to obtain the
properties needed below.

\begin{proposition}[Properties of the Dirichlet solutions]\label{prop:dirichlet-input}
Under the hypotheses of Theorem~\ref{thm:torsion-main}, problem
\eqref{eq:DT} has a unique admissible solution
\[
 U\in C^{3,\gamma}(\overline\Omega)\cap C^\infty(\Omega),\qquad U<0,
\]
with
\begin{equation}\label{eq:torsion-dir-input}
 q_T:=U_\nu>0\ \text{on }\partial\Omega,
 \qquad T_1(D^2U)\ge\vartheta_TI\ \text{on }\overline\Omega,
 \qquad D^2[-\sqrt{-U}]>0\ \text{in }\Omega.
\end{equation}
Assume in addition \eqref{eq:Eg}.  Then
\begin{equation}\label{eq:spectral-dir-input}
 q_E:=(\varphi_D)_\nu>0\ \text{on }\partial\Omega,
 \qquad \varphi_D\in C^{3,\gamma}(\overline\Omega),
 \qquad T_1(D^2\varphi_D)\ge\vartheta_EI\ \text{on }\overline\Omega,
\end{equation}
for some $\vartheta_E>0$, and
\begin{equation}\label{eq:log-dir-input}
 D^2[-\log(-\varphi_D)]>0\qquad\text{in }\Omega.
\end{equation}
Consequently each transformed Hessian has a positive lower bound on
every compact subset of $\Omega$.
\end{proposition}

\begin{proof}
We prove the torsion regularity and convexity first.  For the spectral
problem, regularity and ellipticity have already been obtained in
Lemma~\ref{lem:dirichlet-spectral}; only convexity remains, for which
we use the strong approximation in Lemma~\ref{lem:spectral-approx}.

\smallskip
\noindent\emph{Torsion.}
Choose smooth uniformly convex domains $\Omega_j$ converging to
$\Omega$ in $C^{3,\eta}$ for every $\eta<\gamma$, with uniformly
bounded $C^{3,\gamma}$ charts and a uniform positive lower bound for the
principal curvatures.  On each $\Omega_j$, the classical Dirichlet theory
\cite{CNS,Wang09} provides a smooth admissible solution $U_j$.  Because the right-hand side is the constant $1$, Lyu's global estimate
\cite[Theorem~1.2]{Lyu26} applies uniformly to this family.  Pulling the solutions back to a fixed domain gives uniform
$C^{2,\gamma}$ bounds and, after passage to a subsequence, convergence
in $C^{2,\eta}$.
The limit solves \eqref{eq:DT}; admissibility and uniqueness follow
from the comparison principle.

Since $D^2U\in\overline{\Gamma_2}$ on $\overline\Omega$, we have
\[
\Delta U=\sigma_1(D^2U)\ge0.
\]
On the other hand, by continuity the equation holds on
$\overline\Omega$, and hence
\[
2
=2\sigma_2(D^2U)
=(\Delta U)^2-|D^2U|^2.
\]
Thus $\Delta U$ cannot vanish at any point of $\overline\Omega$.
Consequently
\[
\sigma_1(D^2U)>0,
\qquad
\sigma_2(D^2U)=1>0
\quad\text{on }\overline\Omega,
\]
so
\[
D^2U(x)\in\Gamma_2
\qquad\text{for every }x\in\overline\Omega.
\]
Since $D^2U$ is continuous and $\overline\Omega$ is compact,
$D^2U(\overline\Omega)$ is a compact subset of $\Gamma_2$.
Therefore there exists $\vartheta_T>0$ such that
\[
T_1(D^2U)\ge \vartheta_T I
\qquad\text{on }\overline\Omega.
\]
The equation is now uniformly elliptic in boundary graph coordinates.
Tangential difference quotients satisfy linear Dirichlet equations
with $C^{0,\gamma}$ coefficients.  Boundary Schauder estimates give one additional tangential
derivative.  The equation can then be solved for the normal-normal
second derivative, since the coefficient of that entry is
$T_1(D^2U)[\nu,\nu]>0$.  Hence
$U\in C^{3,\gamma}(\overline\Omega)$.  The Hopf boundary lemma implies
$q_T=U_\nu>0$ on $\partial\Omega$; we set
\[
 q_{T,*}:=\min_{\partial\Omega}q_T>0.
\]

On each smooth approximating domain, $U_j$ is admissible and,
by the positive-density regularity theory and elliptic bootstrapping,
smooth up to the boundary.  The known convexity results therefore give
$D^2[-\sqrt{-U_j}]\ge0$: for $n=3$, use Salani
\cite[Theorem~1.2]{Salani} or Ma--Xu \cite[Theorem~2]{MX}; for
$n\ge4$, use \cite[Theorem~1.3]{LMQS}.  Passing to the limit on compact subsets yields
\[
 D^2[-\sqrt{-U}]\ge0.
\]
For $n=3$, Ma--Xu's local constant-rank lemma
\cite[Lemma~1]{MX} applies because the solution is locally $C^4$ and
admissible, and its transformed Hessian is positive semidefinite.  For $n\ge4$
we use \cite[Theorem~6.2]{LMQS}.  Lemma~\ref{lem:full-rank}(a) excludes
a rank defect and proves the first strict-convexity assertion.

\smallskip
\noindent\emph{Eigenfunction.}
Lemma~\ref{lem:dirichlet-spectral} gives
\eqref{eq:spectral-dir-input} and interior smoothness independently
of any convexity assertion.  It remains to prove
\eqref{eq:log-dir-input}.
Choose the smooth approximating eigenpairs from
Lemma~\ref{lem:spectral-approx}.  For the smooth-domain eigenpairs, the known convexity results give
$D^2[-\log(-\varphi_j)]\ge0$.  In dimension three this follows from
\cite[Theorem~1.1]{LMX} (or Salani's convexity theorem), and for
$n\ge4$ from \cite[Theorem~1.1]{LMQS}.  Because multiplication
by a positive normalization constant does not change the Hessian of the
negative logarithm, passing to the limit in \eqref{eq:spectral-approx} yields
$D^2[-\log(-\varphi_D)]\ge0$ on every compact subset of $\Omega$.
Write $v=-\log(-\varphi_D)$.  Its transformed admissibility is
explicit: since $\varphi_D=-e^{-v}$,
\[
 D^2v-Dv\otimes Dv=e^vD^2\varphi_D\in\Gamma_2
 \qquad\text{in }\Omega.
\]
In dimension three, Liu--Ma--Xu \cite[Theorem~2.1]{LMX} applies to
the limiting transformed equation and implies constant rank: $v$ is
locally $C^4$, transformed-admissible, and has a positive semidefinite
Hessian.  For $n\ge4$ we use \cite[Theorem~4.1]{LMQS}.  Since
$-\log(-\varphi_D)\to+\infty$ at the boundary,
Lemma~\ref{lem:full-rank}(b) forces full rank.  The lower bounds on
compact subsets then follow by continuity.
\end{proof}

\subsection{A parameter-uniform H\"older Robin estimate}
The next estimate will be used for both nonlinear problems.
Proposition~\ref{prop:dirichlet-input} gives
$q_T\in C^{2,\gamma}(\partial\Omega)$ and, under \eqref{eq:Eg},
$q_E\in C^{2,\gamma}(\partial\Omega)$.  Thus the boundary data
$-q_T$ and $-q_E$ in the perturbation equations have exactly the
regularity needed to obtain the $\alpha^{-1}$ factor below.  The extra derivative matters: a uniform $\alpha^{-1}$ gain does not
hold for arbitrary $C^{1,\gamma}$ boundary data.

\begin{theorem}[Uniform Robin resolvent]\label{thm:linear-robin}
Let $n\ge3$, $0<\gamma<1$, and let
$\Omega\subset\mathbb R^n$ be a bounded $C^{3,\gamma}$ domain.  Suppose
\[
 Lw=a^{ij}(x)w_{ij},\qquad a=a^T\in C^{0,\gamma}(\overline\Omega),
 \qquad \vartheta I\le a\le\Theta I.
\]
There exist $\alpha_0\ge1$ and
$C=C(n,\gamma,\Omega,\vartheta,\Theta,\|a\|_{C^{0,\gamma}})>0$
such that, for every $\alpha\ge\alpha_0$, every
$f\in C^{0,\gamma}(\overline\Omega)$,
$g\in C^{2,\gamma}(\partial\Omega)$, and every
$w\in C^{2,\gamma}(\overline\Omega)$ solving
\begin{equation}\label{eq:linear-robin}
 Lw=f\quad\text{in }\Omega,
 \qquad w_\nu+\alpha w=g\quad\text{on }\partial\Omega
\end{equation}
satisfies
\begin{equation}\label{eq:linear-est}
 \|w\|_{C^{2,\gamma}(\overline\Omega)}
 \le C\Bigl(\|f\|_{C^{0,\gamma}(\overline\Omega)}
       +\alpha^{-1}\|g\|_{C^{2,\gamma}(\partial\Omega)}\Bigr).
\end{equation}
For each fixed $\alpha>0$, the problem is uniquely solvable when
$f\in C^{0,\gamma}(\overline\Omega)$ and
$g\in C^{1,\gamma}(\partial\Omega)$.
\end{theorem}

\begin{proof}
We first prove the a priori estimate.  The dependence on $\alpha$ enters
only through the boundary operator, and all localization scales are chosen
independently of $\alpha$.

\smallskip
\noindent\emph{Step 1: the frozen half-space estimate.}
Use coordinates $y\in\mathbb R^{n-1}$ and $t>0$, with outward
normal derivative $-\partial_t$.  For a constant positive definite matrix
$A=(a_{ij})$ set
\[
 L_A=a_{nn}\partial_t^2+2a_{na}\partial_t\partial_{y_a}
                  +a_{ab}\partial_{y_a}\partial_{y_b}.
\]
For $\xi\ne0$ define
\begin{equation}\label{eq:rho}
 \rho(\xi)=
 \frac{\sqrt{a_{nn}a_{ab}\xi_a\xi_b-(a_{na}\xi_a)^2}
        +i\,a_{na}\xi_a}{a_{nn}}.
\end{equation}
Uniform ellipticity implies
$\Re\rho(\xi)\ge c|\xi|$ and $|\rho(\xi)|\le C|\xi|$.
The bounded solution of the homogeneous problem with Robin datum $g$
has tangential Fourier transform
\begin{equation}\label{eq:half-space-multiplier}
 \widehat w(\xi,t)=
 \frac{e^{-t\rho(\xi)}}{\alpha+\rho(\xi)}\,\widehat g(\xi).
\end{equation}
The choice of outward derivative $-\partial_t$ gives the denominator
$\alpha+\rho(\xi)$, whose real part is positive.

Let $g_j$ be a Littlewood--Paley piece supported where
$|\xi|\asymp N=2^j$.  Insert a slightly larger smooth annular cutoff
$\chi$ and rescale $\xi=N\zeta$.  Since $\rho$ is positively
homogeneous of degree one, for every multi-index $\mu$ there exist
$c>0$ and $C_\mu$ such that
\[
 \sup_{|\zeta|\asymp1}
 \left|\partial_\zeta^\mu\left[
 \chi(\zeta)\frac{\alpha+N}{\alpha+N\rho(\zeta)}
           e^{-Nt\rho(\zeta)}\right]\right|
 \le C_\mu e^{-cNt},\qquad t\ge0.
\]
Indeed, $\Re\rho(\zeta)\ge c_0>0$ on the annulus; derivatives of the
rational factor are bounded uniformly in $\alpha/N$, while powers of
$Nt$ arising from derivatives of the exponential are absorbed into a
slightly weaker exponential.  The same estimate holds after inserting
fixed powers of $\rho(\zeta)$ and $\zeta$ for derivatives in $t$ and
$y$.  Rescaling and integrating by parts give a uniform $L^1$
bound for the corresponding convolution kernels, and hence
\begin{equation}\label{eq:dyadic-kernel}
 \|D_y^\beta\partial_t^m w_j(\cdot,t)\|_\infty
 \le C_{m,\beta}\frac{N^{m+|\beta|}}{\alpha+N}
         e^{-cNt}\|g_j\|_\infty ,
 \qquad m\in\mathbb N_0,\ \beta\in\mathbb N_0^{n-1}.
\end{equation}
The constants are uniform with respect to $\alpha/N$.  If
$g\in C^{2,\gamma}$, the dyadic characterization of H\"older spaces
\cite{BCD} gives
$\|g_j\|_\infty\le C N^{-2-\gamma}\|g\|_{C^{2,\gamma}}$.
Set $V_j=D^2w_j$.  Applying \eqref{eq:dyadic-kernel} with one
additional derivative at points $X,Y$ in the closed half-space gives
\begin{equation}\label{eq:dyadic-difference}
 |V_j(X)-V_j(Y)|
 \le C\|g\|_{C^{2,\gamma}}
       \frac{N^{-\gamma}}{\alpha+N}
       \min\{1,N|X-Y|\}.
\end{equation}
Since $(\alpha+N)^{-1}\le\alpha^{-1}$, we sum the terms with
$N\le |X-Y|^{-1}$ using the second factor in the minimum, and the
higher frequencies using the first.  This gives
\[
 [D^2w]_{C^\gamma(\overline{\mathbb R^n_+})}
 \le C\alpha^{-1}\|g\|_{C^{2,\gamma}}.
\]
Summation without taking differences also controls the $C^0$ norm
of $D^2w$ and the lower derivatives.  For low frequencies, let $\rho(D)$ denote the tangential Fourier
multiplier with symbol $\rho(\xi)$ and use the representation
\[
 (\alpha+\rho(D))^{-1}
   =\int_0^\infty e^{-\alpha s}e^{-s\rho(D)}\,ds.
\]
In fact,
\[
 \rho(\xi)=\sqrt{\xi^TQ\xi}+i b_0\cdot\xi,\qquad
 Q_{ab}=\frac{a_{ab}}{a_{nn}}-
                \frac{a_{na}a_{nb}}{a_{nn}^2},\qquad
 (b_0)_a=\frac{a_{na}}{a_{nn}}.
\]
The matrix $Q$ is positive definite with bounds determined by
$\vartheta,\Theta$.  Thus $e^{-s\rho(D)}$ is an anisotropic Poisson
operator followed by translation.  Its kernel has mass one, so its
$L^\infty$ norm is at most one.  A fixed low-frequency cutoff and its
derivatives have integrable kernels; the same holds after
multiplication by the powers of $\rho$ needed for normal derivatives.
Consequently the low-frequency terms, including their second-derivative
H\"older seminorms, are bounded by
$C\alpha^{-1}\|g\|_{C^{2,\gamma}}$.  Together with the annular
estimates, this proves
\begin{equation}\label{eq:frozen-boundary-est}
 \|w\|_{C^{2,\gamma}(\overline{\mathbb R^n_+})}
 \le C\alpha^{-1}\|g\|_{C^{2,\gamma}(\mathbb R^{n-1})}
 \qquad(\alpha\ge1).
\end{equation}

For the localized interior source $F$, fix $R_0$ so that
$\operatorname{supp}F\subset\overline{B_{R_0}^+}$.  Let $z$ be the
zero-Dirichlet solution of $L_Az=F$ in the half-space.  A linear change
of coordinates reduces its Green kernel to the reflected Laplace
kernel.  The bound by $C|X-Y|^{2-n}$ gives
$\|z\|_\infty\le CR_0^2\|F\|_\infty$.  Boundary and interior
Schauder estimates on unit half-balls and balls then yield
\begin{equation}\label{eq:frozen-source}
 \|z\|_{C^{2,\gamma}}+\|z_\nu\|_{C^{1,\gamma}}
 \le C_{R_0}\|F\|_{C^{0,\gamma}}.
\end{equation}
Here $C_{R_0}$ depends only on
$n,\gamma,\vartheta,\Theta,R_0$, and not on $\alpha$.

For clarity, let $\mathcal P_\alpha h$ denote the homogeneous
solution with Robin datum $h$.  If $h\in C^{1,\gamma}$, then
$\|h_j\|_\infty\le C N^{-1-\gamma}\|h\|_{C^{1,\gamma}}$.
Equation~\eqref{eq:dyadic-kernel} gives
\begin{align*}
 &|D^2\mathcal P_\alpha h_j(X)
                 -D^2\mathcal P_\alpha h_j(Y)|\\
 &\qquad\le C\|h\|_{C^{1,\gamma}}
       \frac{N^{1-\gamma}}{\alpha+N}
                     \min\{1,N|X-Y|\}\\
 &\qquad\le C\|h\|_{C^{1,\gamma}}N^{-\gamma}
                     \min\{1,N|X-Y|\}.
\end{align*}
Summing below and above $N=|X-Y|^{-1}$, and using the preceding
Poisson representation for low frequencies, proves
\begin{equation}\label{eq:frozen-rough-data}
 \|\mathcal P_\alpha h\|_{C^{2,\gamma}}
       \le C\|h\|_{C^{1,\gamma}}\qquad(\alpha\ge1).
\end{equation}
Since the solution with source $F$ and datum $g$ is
$z+\mathcal P_\alpha(g-z_\nu)$, estimates
\eqref{eq:frozen-boundary-est}, \eqref{eq:frozen-source}, and
\eqref{eq:frozen-rough-data} give
\begin{equation}\label{eq:frozen-combined}
 \|w\|_{C^{2,\gamma}}
 \le C_{R_0}\bigl(\|F\|_{C^{0,\gamma}}
                      +\alpha^{-1}\|g\|_{C^{2,\gamma}}\bigr).
\end{equation}
The source term is bounded uniformly, whereas the prescribed smoother
boundary datum carries the additional factor $\alpha^{-1}$.

\smallskip
\noindent\emph{Step 2: localization.}
Near the boundary use normal coordinates
\[
 x=\Psi(y,t)=X(y)-t\nu(X(y)).
\]
The $C^{3,\gamma}$ regularity of $\partial\Omega$ gives
$C^{2,\gamma}$ normal coordinates on a fixed boundary neighborhood.
The principal coefficients in these coordinates have bounded
$C^{0,\gamma}$ norms, as do the lower-order coefficients.  At $t=0$ the boundary condition is exactly
$-w_t+\alpha w=g$.

Choose cutoffs of the form
$\chi(y,t)=\chi_T(y)\chi_N(t)$ with $\chi_N'(0)=0$.
The cutoff introduces no additional boundary term:
\[
 (-\partial_t+\alpha)(\chi w)=\chi g\quad\text{on }t=0.
\]
Freeze the principal matrix at the center of a sufficiently small
chart and apply \eqref{eq:frozen-combined} to the localized function.
To track the lower-order terms, write the transformed operator as
$\widetilde L=\widetilde a^{ij}\partial_{ij}
+\widetilde b^i\partial_i$, let $\widetilde w=w\circ\Psi$, and put
$v=\chi\widetilde w$.  With $A_0$ the frozen principal matrix,
\begin{equation}\label{eq:localized-equation}
 L_{A_0}v=\chi\widetilde f
 +(A_0-\widetilde a):D^2v-\widetilde b\cdot Dv
                         +[\widetilde L,\chi]\widetilde w.
\end{equation}
The commutator equals
$2\widetilde a^{ij}\chi_i\widetilde w_j+
(\widetilde L\chi)\widetilde w$.  It contains at most first
 derivatives of $\widetilde w$, and its $C^{0,\gamma}$ norm is at
most $C_r\|\widetilde w\|_{C^{1,\gamma}}$ on a chart of fixed
radius $r$.  No commutator term contains $\alpha$.
Dropping tildes in the principal coefficient estimate, we have
\[
 \|(a-a_0)D^2(\chi w)\|_{C^{0,\gamma}}
 \le C\operatorname{osc}(a)\|\chi w\|_{C^{2,\gamma}}
    +C[a]_{C^\gamma}\|D^2(\chi w)\|_{C^0},
\]
where $\operatorname{osc}(a)$ is the oscillation of the coefficient
matrix on the chart and $[a]_{C^\gamma}$ is its H\"older seminorm.
Choose the chart radius so that the contribution of the
coefficient oscillation to $\|\chi w\|_{C^{2,\gamma}}$ has coefficient
at most $1/4$.  With this radius fixed, use
\[
 \|D^2v\|_{C^0}+\|v\|_{C^{1,\gamma}}
 \le\varepsilon\|v\|_{C^{2,\gamma}}
      +C_\varepsilon\|v\|_{C^0}
\]
and choose $\varepsilon$ small enough to absorb the remaining
$C^{2,\gamma}$ term into the left-hand side.  After summing over a fixed finite boundary cover, the commutator
and coordinate terms are bounded by $C_r\|w\|_{C^{1,\gamma}}$.
A second application of interpolation absorbs this term, up to
$C\|w\|_{C^0}$, into the global $C^{2,\gamma}$ norm.  The radii,
cutoffs, and both interpolation parameters are fixed before $\alpha$
is chosen.  Together with the interior Schauder estimates, this gives
\begin{equation}\label{eq:localized-robin}
 \|w\|_{C^{2,\gamma}}
 \le C\Bigl(\|f\|_{C^{0,\gamma}}
    +\alpha^{-1}\|g\|_{C^{2,\gamma}}
    +\|w\|_{C^0}\Bigr),
\end{equation}
where $C$ is independent of all sufficiently large $\alpha$.

\smallskip
\noindent\emph{Step 3: the zeroth-order term.}
Let $p$ solve $Lp=-1$ in $\Omega$ with $p=0$ on
$\partial\Omega$.  Then $p\ge0$, and its $C^{2,\gamma}$ norm is bounded
in terms of the domain and the coefficients of $L$.  Put
$F=\|f\|_\infty$, $G=\|g\|_\infty$, and
\[
 W=Fp+\frac{G+F\|p_\nu\|_\infty}{\alpha}.
\]
We have $LW=-F$ and
$W_\nu+\alpha W\ge G$.  Applying the maximum principle to both $w-W$ and $-w-W$ gives
\begin{equation}\label{eq:C0-robin}
 \|w\|_\infty\le CF+\alpha^{-1}G.
\end{equation}
Combining this with \eqref{eq:localized-robin} proves
\eqref{eq:linear-est}.

For each fixed $\alpha>0$, classical oblique Schauder theory
\cite{Lieberman} makes the boundary-value operator Fredholm.  This
solvability argument is independent of the large-$\alpha$ estimate.  A homotopy from the Laplacian preserves obliqueness and positivity
of the Robin coefficient, and therefore preserves Fredholm index zero.  The maximum principle implies that the kernel is trivial;
hence the operator is bijective and the solution is unique.
\end{proof}

\begin{remark}[The boundary regularity behind the $\alpha^{-1}$ gain]
\label{rem:robin-symbol}
The frozen half-space symbol explains the regularity requirement
on the Robin datum in Theorem~\ref{thm:linear-robin}.  For a dyadic piece with $|\xi|\asymp N$,
\eqref{eq:half-space-multiplier} gives the schematic bound
\[
 \|D^2w_j\|_\infty
 \lesssim \frac{N^2}{\alpha+N}\,\|g_j\|_\infty.
\]
If $g\in C^{2,\gamma}$, then
$\|g_j\|_\infty\lesssim N^{-2-\gamma}\|g\|_{C^{2,\gamma}}$;
after multiplication by the H\"older weight $N^\gamma$, the dyadic
contribution is bounded by $(\alpha+N)^{-1}$ and hence by
$\alpha^{-1}$.  This yields the factor appearing in \eqref{eq:linear-est}.  For general $g\in C^{1,\gamma}$, the weighted factor is instead
$N/(\alpha+N)\le1$.  This gives a uniform estimate without a uniform
$\alpha^{-1}$ improvement.  Equivalently,
\[
 g\in C^{2,\gamma}\quad\Longrightarrow\quad O(\alpha^{-1}),
 \qquad
 g\in C^{1,\gamma}\quad\Longrightarrow\quad O(1)
\]
at the level of the second-derivative H\"older norm.  The failure of
an $\alpha^{-1}$ estimate for general $C^{1,\gamma}$ data is already
visible for $L=\Delta$.  Take $\alpha=N$ and
\[
 h_N(y)=N^{-1-\gamma}\cos(Ny_1),\qquad
 w_N(y,t)=\frac{N^{-1-\gamma}}{2N}e^{-Nt}\cos(Ny_1).
\]
Then $\|h_N\|_{C^{1,\gamma}}\le C$, while
$[D^2w_N]_{C^\gamma}\ge c>0$ independently of $N$.
Thus no additional $\alpha^{-1}$ factor is possible at this lower
data regularity, even in the frozen half-space model.
\end{remark}

\subsection{The augmented spectral inverse}
The Dirichlet spectral operators $L_E$, $M_E$, and $\mathfrak n$ were
defined in \eqref{eq:linearizations}.  Their zero-Dirichlet augmented
inverse is already available from Lemma~\ref{lem:dirichlet-spectral}.
We now pass from Dirichlet to Robin boundary conditions.  For torsion,
set
\[
 L_Tv=T_1(D^2U):D^2v.
\]
Proposition~\ref{prop:dirichlet-input} ensures that both $L_T$ and
$L_E$ satisfy the hypotheses of Theorem~\ref{thm:linear-robin}.

\begin{lemma}[Augmented spectral resolvent]\label{lem:augmented}
Assume \eqref{eq:Eg}.  There exist $\alpha_1,C>0$ such that,
for every $\alpha\ge\alpha_1$ and every
$f\in C^{0,\gamma}(\overline\Omega)$,
$g\in C^{2,\gamma}(\partial\Omega)$, and $s\in\mathbb R$,
the system
\begin{equation}\label{eq:augmented}
 M_Ev-\ell\varphi_D^2=f,
 \qquad v_\nu+\alpha v=g,
 \qquad \mathfrak n(v)=s
\end{equation}
has a unique solution $(v,\ell)\in
C^{2,\gamma}(\overline\Omega)\times\mathbb R$ and
\begin{equation}\label{eq:augmented-est}
 \|v\|_{C^{2,\gamma}(\overline\Omega)}+|\ell|
 \le C\Bigl(\|f\|_{C^{0,\gamma}(\overline\Omega)}+|s|
      +\alpha^{-1}\|g\|_{C^{2,\gamma}(\partial\Omega)}\Bigr).
\end{equation}
\end{lemma}

\begin{proof}
Write $\phi=\varphi_D$ and $a_E=T_1(D^2\phi)$.  By
Lemma~\ref{lem:dirichlet-spectral}, $a_E\in C^{1,\gamma}$ and the
map \eqref{eq:dir-aug-map} is an isomorphism for $\eta=\gamma$.
It remains to perturb this inverse to the Robin boundary condition.

Let $R_D$ and $R_\alpha$ denote the inverses of $L_E$ with,
respectively, zero Dirichlet and homogeneous Robin boundary conditions.  If $h\in C^{2,\gamma}$, then $a_E\in C^{1,\gamma}$ and Dirichlet
Schauder theory give $R_Dh\in C^{3,\gamma}$.  Consequently, $z=(R_\alpha-R_D)h$ is $L_E$-harmonic, with Robin
datum $-\partial_\nu R_Dh\in C^{2,\gamma}$.  Theorem~\ref{thm:linear-robin} therefore gives
\begin{equation}\label{eq:R-convergence}
 \|R_\alpha-R_D\|_{\mathcal L(C^{2,\gamma},C^{2,\gamma})}
 \le C\alpha^{-1}.
\end{equation}

On $\mathbb X=C^{2,\gamma}(\overline\Omega)\times\mathbb R$ define
\[
 \mathcal A_\alpha(v,\ell)
 =\left(v-R_\alpha(2\lambda_D\phi v+\ell\phi^2),
        \mathfrak n(v)\right).
\]
Define $\mathcal A_D$ by replacing $R_\alpha$ with $R_D$.  The first component of $\mathcal A_D$ need not have zero boundary
trace.  We compare its invertibility with that of the zero-Dirichlet
augmented map.  Given
$(h,s)\in\mathbb X$, the equation
$\mathcal A_D(v,\ell)=(h,s)$ is equivalent, after writing $v=h+u$, to
\begin{equation}\label{eq:AD-equivalence}
 \begin{cases}
 M_Eu-\ell\phi^2=2\lambda_D\phi h&\text{in }\Omega,\\
 u=0&\text{on }\partial\Omega,\\
 \mathfrak n(u)=s-\mathfrak n(h).
 \end{cases}
\end{equation}
Indeed, the first component of $\mathcal A_D(v,\ell)=(h,s)$ is
equivalent to $u=R_D(2\lambda_D\phi(h+u)+\ell\phi^2)$.  Applying $L_E$
yields \eqref{eq:AD-equivalence}; applying $R_D$ proves the converse.
The zero-Dirichlet isomorphism \eqref{eq:dir-aug-map} determines a
unique $(u,\ell)$ and hence $v=h+u$.  The resulting inverse is bounded
on $\mathbb X$.  Thus $\mathcal A_D$ is invertible.  Because
multiplication by $\phi$ is bounded on $C^{2,\gamma}$,
\eqref{eq:R-convergence} implies
\[
 \|\mathcal A_\alpha-\mathcal A_D\|_{\mathcal L(\mathbb X)}
 \le C\alpha^{-1}.
\]
For sufficiently large $\alpha$, the Neumann series therefore
inverts $\mathcal A_\alpha$ with a uniform bound.  Finally, Theorem~\ref{thm:linear-robin} bounds $R_\alpha f$ by
$C\|f\|_{C^{0,\gamma}}$ and the Robin lifting of $g$ by
$C\alpha^{-1}\|g\|_{C^{2,\gamma}}$.  Applying
$\mathcal A_\alpha^{-1}$ proves \eqref{eq:augmented-est}.
\end{proof}

\section{Large-parameter Robin solutions}\label{sec:branches}

We work in the spaces
\[
 X=C^{2,\gamma}(\overline\Omega),\qquad
 Y=C^{0,\gamma}(\overline\Omega).
\]
We equip $X\times\mathbb R$ with the norm
\[
 \|(v,\mu)\|_{X\times\mathbb R}=\|v\|_X+|\mu|.
\]
Fix $C_{\rm emb}\ge1$ such that, for every $v\in X$,
\[
 \|v\|_{C^0}+\|Dv\|_{C^0}+\|D^2v\|_{C^0}
 \le C_{\rm emb}\|v\|_X,
\]
where the matrix norm on $D^2v$ is the Frobenius norm.  We also use
the fact that $Y$ is an algebra.

\subsection{Exact difference equations and quadratic estimates}

For the torsion equation, set
\[
 v=u_\alpha-U.
\]
By the exact quadratic identity \eqref{eq:quad},
\begin{equation}\label{eq:torsion-diff}
 L_Tv=-\sigma_2(D^2v)\quad\text{in }\Omega,
 \qquad v_\nu+\alpha v=-q_T\quad\text{on }\partial\Omega.
\end{equation}
For the eigenvalue problem set
\[
 \phi=\varphi_D,\qquad \lambda=\lambda_D,\qquad
 e=\varphi_\alpha-\phi,\qquad
 \mu=\lambda_\alpha-\lambda .
\]
Expanding both sides of the eigenvalue equation gives
\begin{equation}\label{eq:eigen-diff}
 M_Ee-\mu\phi^2
 =-\sigma_2(D^2e)+\lambda e^2+2\mu\phi e+\mu e^2
 =:\mathcal R(e,\mu),
\end{equation}
with
\begin{equation}\label{eq:eigen-diff-boundary}
 e_\nu+\alpha e=-q_E.
\end{equation}
Expanding the normalization gives
\begin{align}
 1
 &=\int_\Omega(-\phi-e)^3
  =1-3\int_\Omega\phi^2e
    -3\int_\Omega\phi e^2-\int_\Omega e^3,
\end{align}
hence
\begin{equation}\label{eq:norm-rem}
 \mathfrak n(e)
 =-\int_\Omega\phi e^2-\frac13\int_\Omega e^3
 =:\mathcal N(e).
\end{equation}

For some constant $C_X$ depending only on the fixed Dirichlet data,
\begin{align}
 \|\sigma_2(D^2v)\|_Y
 &\le C_X\|v\|_X^2,\label{eq:quad-est}\\
 \|\sigma_2(D^2v)-\sigma_2(D^2\widetilde v)\|_Y
 &\le C_X(\|v\|_X+\|\widetilde v\|_X)
          \|v-\widetilde v\|_X.\label{eq:quad-lip}
\end{align}
Likewise, on the unit ball of $X\times\mathbb R$,
\begin{equation}\label{eq:spectral-quad}
 \|\mathcal R(e,\mu)\|_Y+|\mathcal N(e)|
 \le C_X(\|e\|_X+|\mu|)^2,
\end{equation}
and
\begin{align}
 &\|\mathcal R(e,\mu)-\mathcal R(\widetilde e,\widetilde\mu)\|_Y
 +|\mathcal N(e)-\mathcal N(\widetilde e)|\notag\\
 &\qquad\le
 C_X\bigl(\|e\|_X+|\mu|+\|\widetilde e\|_X+|\widetilde\mu|\bigr)
 \bigl(\|e-\widetilde e\|_X+|\mu-\widetilde\mu|\bigr).
 \label{eq:spectral-lip}
\end{align}
The H\"older product estimates give these bounds:
$\sigma_2(D^2v)$ is quadratic in the entries of $D^2v$, and every term
in $\mathcal R$ and $\mathcal N$ has degree at least two in $(e,\mu)$.
Thus the interior nonlinearities are quadratic, whereas the Robin
boundary data $-q_T$ and $-q_E$ are of order one.  By
Proposition~\ref{prop:dirichlet-input}, these data belong to
$C^{2,\gamma}(\partial\Omega)$, precisely the class for which the
linear Robin estimates give an $O(\alpha^{-1})$ contribution.  We therefore seek fixed points in balls of radius
$R/\alpha$.

\subsection{The contraction argument}

For torsion, define $\mathcal T_\alpha:X\to X$ through
\begin{equation}\label{eq:Ttorsion-def}
 \begin{cases}
 L_T(\mathcal T_\alpha v)=-\sigma_2(D^2v)&\text{in }\Omega,\\
 (\mathcal T_\alpha v)_\nu+\alpha\mathcal T_\alpha v=-q_T
   &\text{on }\partial\Omega.
 \end{cases}
\end{equation}
Theorem~\ref{thm:linear-robin} and
\eqref{eq:quad-est}--\eqref{eq:quad-lip} imply
\begin{align}
 \|\mathcal T_\alpha v\|_X
 &\le C_0\alpha^{-1}+C_1\|v\|_X^2,\label{eq:Ttorsion-self}\\
 \|\mathcal T_\alpha v-\mathcal T_\alpha\widetilde v\|_X
 &\le C_2(\|v\|_X+\|\widetilde v\|_X)
            \|v-\widetilde v\|_X.\label{eq:Ttorsion-lip}
\end{align}

For the eigenvalue problem, apply Lemma~\ref{lem:augmented}.
Given $(e,\mu)$, define
$\mathcal T_\alpha^E(e,\mu)=(\widehat e,\widehat\mu)$ as the solution of
\begin{equation}\label{eq:Teigen-def}
 \begin{cases}
 M_E\widehat e-\widehat\mu\,\phi^2=\mathcal R(e,\mu),\\
 \widehat e_\nu+\alpha\widehat e=-q_E,\\
 \mathfrak n(\widehat e)=\mathcal N(e).
 \end{cases}
\end{equation}
The estimates \eqref{eq:augmented-est},
\eqref{eq:spectral-quad}, and \eqref{eq:spectral-lip} give
\begin{align}
 \|\mathcal T_\alpha^E(e,\mu)\|_{X\times\mathbb R}
 &\le C_0\alpha^{-1}
   +C_1(\|e\|_X+|\mu|)^2,\label{eq:Teigen-self}\\
 \|\mathcal T_\alpha^E(e,\mu)
      -\mathcal T_\alpha^E(\widetilde e,\widetilde\mu)\|_{X\times\mathbb R}
 &\le C_2\bigl(\|e\|_X+|\mu|
       +\|\widetilde e\|_X+|\widetilde\mu|\bigr)\notag\\
 &\qquad\times
 \bigl(\|e-\widetilde e\|_X+|\mu-\widetilde\mu|\bigr).
 \label{eq:Teigen-lip}
\end{align}

Choose
\[
 R=\max\{1,2C_0\}.
\]
For torsion, take $\alpha$ above the threshold in
Theorem~\ref{thm:linear-robin}; for the spectral problem, also impose
the threshold in Lemma~\ref{lem:augmented}.  The constants above,
and hence $R$, may be chosen separately for the two problems.
In each case require
\begin{equation}\label{eq:contraction-threshold}
 \alpha\ge\max\{R,\,2C_1R,\,4C_2R\}.
\end{equation}
If $\|v\|_X\le R/\alpha$, then
\[
 C_0\alpha^{-1}+C_1R^2\alpha^{-2}
 \le R\alpha^{-1}.
\]
The same self-mapping estimate holds for the spectral map on
$X\times\mathbb R$.  By \eqref{eq:Ttorsion-lip} and
\eqref{eq:Teigen-lip}, both Lipschitz constants are at most $1/2$.
The contraction mapping theorem gives a unique fixed point for each
map in its corresponding ball.
In particular,
\begin{equation}\label{eq:branch-rates}
 \|u_\alpha-U\|_{C^{2,\gamma}}
 \le R\alpha^{-1},
 \qquad
 \|\varphi_\alpha-\phi\|_{C^{2,\gamma}}
   +|\lambda_\alpha-\lambda|
 \le R\alpha^{-1}.
\end{equation}
These fixed points solve the perturbation equations in a
$C^{2,\gamma}$ neighborhood of the Dirichlet solutions.  We must still
verify negativity and admissibility, and prove uniqueness among all
negative admissible solutions rather than only within the fixed-point
balls.

\subsection{The torsion solutions}

Since $U\in C^{3,\gamma}$ and
$D^2U(\overline\Omega)$ is a compact subset of
$\Gamma_2$, set
\begin{equation}\label{eq:Gamma-gap}
 d_\Gamma=
 \min_{x\in\overline\Omega}
 \dist(D^2U(x),\partial\Gamma_2)>0.
\end{equation}
With the constant $C_{\rm emb}$ fixed above, if
\[
 \alpha\ge \frac{2C_{\rm emb}R}{d_\Gamma},
\]
then the torsion fixed point satisfies
$D^2u_\alpha\in\Gamma_2$ on the closure.  In particular
$\Delta u_\alpha>0$.  Suppose that the maximum $M$ of $u_\alpha$ on $\overline\Omega$
were nonnegative.  Since $\Delta u_\alpha>0$, a maximum occurs on
$\partial\Omega$, where $(u_\alpha)_\nu\ge0$.  If $M>0$, the Robin
condition gives $(u_\alpha)_\nu=-\alpha M<0$, a contradiction.
If $M=0$, the Robin condition gives $(u_\alpha)_\nu=0$, whereas the
Hopf lemma yields $(u_\alpha)_\nu>0$ because the solution is
nonconstant.  Hence
\begin{equation}\label{eq:torsion-negative}
 u_\alpha<0\quad\text{on }\overline\Omega.
\end{equation}

We have shown that $u_\alpha$ is a negative admissible classical
solution of \eqref{eq:RT}, satisfying \eqref{eq:torsion-rate-main}.
For uniqueness, let $u$ and $\widetilde u$ be two admissible classical
solutions with the same $\alpha$.  Their Hessians extend continuously
to $\overline\Omega$, and the identity
$2\sigma_2(D^2u)=(\Delta u)^2-|D^2u|^2=2$, together with
$\Delta u\ge0$, gives $D^2u\in\Gamma_2$ on $\overline\Omega$;
the same holds for $\widetilde u$.  Convexity of $\Gamma_2$ and
compactness of $\overline\Omega$ therefore give a constant
$\theta=\theta(u,\widetilde u)>0$ such that, for
$w=u-\widetilde u$,
\begin{equation}\label{eq:torsion-unique-eq}
 a^{ij}w_{ij}=0,\qquad
 a=T_1\!\left(\frac{D^2u+D^2\widetilde u}{2}\right)
 \ge\theta I,\qquad w_\nu+\alpha w=0.
\end{equation}
The maximum principle applied to $w$ and $-w$ gives $w=0$.
Only the transformed-Hessian conclusion of
Theorem~\ref{thm:torsion-main} remains to be proved.

\subsection{The normalized spectral solutions}

Assume \eqref{eq:Eg} and set
\[
 a_E=T_1(D^2\phi),\qquad
 a_E\ge\vartheta_EI.
\]
Since
$\tr a_E=(n-1)\Delta\phi$, there exists
$\tau_E>0$ such that
\begin{equation}\label{eq:laplace-gap}
 \Delta\phi\ge\tau_E\quad\text{on }\overline\Omega.
\end{equation}
Choose $\delta_0>0$ smaller than both the tubular radius of
$\partial\Omega$ and the inradius of $\Omega$ so that, for
$x=y-d\nu(y)$ with $0\le d\le\delta_0$,
\begin{equation}\label{eq:phi-normal-collar}
 D\phi(x)\cdot\nu(y)\ge\frac12q_{E,*}.
\end{equation}
On the compact interior set
\[
 K_{\delta_0}=\{x\in\overline\Omega:d(x)\ge\delta_0\},\qquad
 m_{\phi,\delta_0}=\min_{x\in K_{\delta_0}}(-\phi(x))>0.
\]

Let $(e,\mu)$ be the spectral fixed point.  Increase the parameter
threshold so that
\begin{equation}\label{eq:spectral-sign-threshold}
 C_{\rm emb}R\alpha^{-1}
 \le \min\left\{\frac14q_{E,*},
                 \frac12m_{\phi,\delta_0},
                 \frac{\lambda_D}{2},
                 \frac{\tau_E}{2n}\right\}.
\end{equation}
Then $\lambda_\alpha=\lambda_D+\mu\ge\lambda_D/2$, and the outward
normal component of $D(\phi+e)$ is positive whenever
$d(x)\le\delta_0$.  The Robin
condition implies
\[
 \varphi_\alpha|_{\partial\Omega}
 =-\frac{(\varphi_\alpha)_\nu}{\alpha}<0.
\]
For $x=y-d\nu(y)$ in this neighborhood,
\[
 \varphi_\alpha(x)=\varphi_\alpha(y)
 -\int_0^dD\varphi_\alpha(y-t\nu(y))\cdot\nu(y)\,dt<0.
\]
On $K_{\delta_0}$, the $C^0$ part of
\eqref{eq:spectral-sign-threshold} implies
$\varphi_\alpha<0$ on $K_{\delta_0}$.  Therefore
\begin{equation}\label{eq:spectral-negative}
 \varphi_\alpha<0\quad\text{on }\overline\Omega.
\end{equation}

The threshold above, together with \eqref{eq:laplace-gap}, gives
$\Delta\varphi_\alpha>0$.  Since the equation at the fixed point is
\[
 \sigma_2(D^2\varphi_\alpha)
 =\lambda_\alpha\varphi_\alpha^2>0,
\]
we have
\[
 D^2\varphi_\alpha\in\Gamma_2
 \quad\text{on }\overline\Omega.
\]
Increasing the lower bound on $\alpha$ if necessary, linearity of
$T_1$ gives the uniform ellipticity bound
\begin{equation}\label{eq:spectral-elliptic-branch}
 T_1(D^2\varphi_\alpha)\ge\frac12\vartheta_E I.
\end{equation}
The third equation in \eqref{eq:Teigen-def} is the normalization
condition.  Hence the fixed point is a normalized eigenpair and
satisfies \eqref{eq:eigen-rate-main}.

The contraction gives uniqueness only in an $O(\alpha^{-1})$
neighborhood of the Dirichlet pair.  Uniqueness in the class specified
in Theorem~\ref{thm:eigen-main} follows from the next comparison
result, which does not require $\alpha$ to be large.

\begin{proposition}[Spectral uniqueness]\label{prop:spectral-unique}
For each fixed $\alpha>0$, at most one positive eigenvalue admits
a negative admissible classical Robin eigenfunction.  Such
eigenfunctions are unique up to multiplication by a positive constant.
\end{proposition}
\begin{proof}
Let $(\lambda_i,\varphi_i)$, $i=1,2$, be two such pairs.  We first
show that $\varphi_i<0$ on $\partial\Omega$ as well.  The interior
negativity and continuity give $\varphi_i\le0$ on $\partial\Omega$.  If $\varphi_i(y)=0$ at a boundary point, then $y$ is a boundary
maximum.  Admissibility gives $\Delta\varphi_i>0$ in $\Omega$, so the
Hopf lemma implies $(\varphi_i)_\nu(y)>0$.  The Robin condition,
however, gives $(\varphi_i)_\nu(y)=0$, a contradiction.

Consequently
\[
 V_i=-\log(-\varphi_i)
\]
belongs to $C^2(\overline\Omega)$ and satisfies $(V_i)_\nu=\alpha$.
Set
\[
 Z=V_1-V_2,\qquad
 B_i=D^2V_i-DV_i\otimes DV_i.
\]
Since $D^2\varphi_i=e^{-V_i}B_i$, homogeneity of $\sigma_2$ implies
\[
 F(B_i)=\sqrt{\lambda_i},\qquad F=\sqrt{\sigma_2}.
\]
We claim that $B_i\in\Gamma_2$ on $\overline\Omega$.  This holds
in the interior by admissibility.  On the boundary, the strict
negativity just proved gives
\[
 \sigma_2(D^2\varphi_i)=\lambda_i\varphi_i^2>0.
\]
Continuity from the interior also gives $\Delta\varphi_i\ge0$ on
the boundary.  If $\Delta\varphi_i=0$ there, the identity
$2\sigma_2(A)=(\operatorname{tr}A)^2-|A|^2$ would imply
$\sigma_2(D^2\varphi_i)\le0$, a contradiction.  Therefore
$D^2\varphi_i\in\Gamma_2$ on $\overline\Omega$.  Since
$B_i=e^{V_i}D^2\varphi_i$ and $e^{V_i}>0$, the same conclusion holds for
$B_i$.

For $V_t=(1-t)V_2+tV_1$, direct calculation gives
\begin{equation}\label{eq:B-path}
 B_t:=D^2V_t-DV_t\otimes DV_t
 =(1-t)B_2+tB_1+t(1-t)DZ\otimes DZ.
\end{equation}
Convexity of the G{\aa}rding cone places the sum of the first two
terms in $\Gamma_2$.  If $A\in\Gamma_2$ and $P\ge0$, then
\[
 \sigma_1(A+P)>0,
 \qquad
 \sigma_2(A+P)
 =\sigma_2(A)+T_1(A):P+\sigma_2(P)>0,
\]
because $T_1(A)>0$ and $\sigma_2(P)\ge0$.  Thus $B_t\in\Gamma_2$ for all
$(x,t)\in\overline\Omega\times[0,1]$.  Continuity of $B_t$ and compactness of
$\overline\Omega\times[0,1]$ then show that the image of this path is
compactly contained in $\Gamma_2$.

The fundamental theorem of calculus along this path yields
\begin{equation}\label{eq:Z-linear}
 A^{ij}Z_{ij}+b^jZ_j
 =\sqrt{\lambda_1}-\sqrt{\lambda_2},
 \qquad Z_\nu=0,
\end{equation}
where $F^{ij}(B)$ denotes the Frobenius derivative on
$\operatorname{Sym}(n)$, characterized for every symmetric $H$ by
$DF(B)[H]=F^{ij}(B)H_{ij}$, and
\[
 A^{ij}=\int_0^1F^{ij}(B_t)\,dt,\qquad
 b^j=-2\int_0^1F^{ij}(B_t)(V_t)_i\,dt.
\]
As $F^{ij}(B)=T_1(B)^{ij}/(2\sqrt{\sigma_2(B)})$ and the path is
compactly contained in $\Gamma_2$, the matrix $A=(A^{ij})$ is
uniformly elliptic and the vector $b$ is bounded.

If the constant on the right of \eqref{eq:Z-linear} were positive,
$Z$ could not be constant.  The maximum principle would then force a maximum to occur on
$\partial\Omega$.  The Hopf lemma gives $Z_\nu>0$ there, contrary to
$Z_\nu=0$.  Applying the argument to $-Z$ excludes a negative constant.  Therefore $\lambda_1=\lambda_2$.  When the right-hand side vanishes, the strong maximum principle and
the Hopf lemma show that $Z$ is constant.  Hence
$\varphi_1=e^{-Z}\varphi_2$, so the two eigenfunctions are positive multiples of one another.
\end{proof}

The normalization in \eqref{eq:RE} fixes this factor.  Hence the
spectral fixed point is the unique normalized negative admissible
eigenpair, with the convergence stated in
Theorem~\ref{thm:eigen-main}.

We have established existence, admissibility, uniqueness, and
$O(\alpha^{-1})$ convergence for both Robin problems.  The remaining
task is to prove positivity of the transformed Hessians near
$\partial\Omega$.

\section{From \texorpdfstring{$C^2$}{C2}-stability to concavity}\label{sec:concavity}

Estimate \eqref{eq:branch-rates} gives $C^2$ convergence, but the
Dirichlet solutions vanish on $\partial\Omega$, where the coefficients
in both transformed Hessians become singular.  We use a common
criterion for the two transforms.

\subsection{A unified Hessian criterion}
For $f<0$ and $c>0$, define
\begin{equation}\label{eq:Kc}
 \mathscr K_c[f]=D^2f+c\frac{Df\otimes Df}{-f}.
\end{equation}
Direct differentiation yields
\begin{equation}\label{eq:transform-formulas}
 D^2[-\sqrt{-f}]
 =\frac1{2\sqrt{-f}}\mathscr K_{1/2}[f],
 \qquad
 D^2[-\log(-f)]
 =\frac1{-f}\mathscr K_1[f].
\end{equation}
Consequently, positivity of $\mathscr K_c[f_\alpha]$ for
$c=1/2$ and $c=1$, respectively, suffices for the two transforms.

\subsection{A Schur-complement lemma}

\begin{lemma}[Uniform Schur positivity]\label{lem:schur}
Let $a,L,R,q_-,q_+$ be positive and $0<c_-\le c_+$.
In the decomposition $\nu^\perp\oplus\mathbb R\nu$, suppose
\[
 D^2f=\begin{pmatrix}A&b\\ b^T&r\end{pmatrix},
 \qquad Df=p+q\nu,\qquad w=-f>0,
\]
where
\[
 A\ge aI,\qquad |b|+|r|\le L,\qquad
 q_-\le q\le q_+,\qquad |p|\le Rw.
\]
There is $w_0>0$, depending only on the stated constants, such that,
whenever $w\le w_0$ and $c\in[c_-,c_+]$,
\begin{equation}\label{eq:K-lower}
 \mathscr K_c[f][\zeta+s\nu,\zeta+s\nu]
 \ge \frac a2|\zeta|^2+\frac{c_-q_-^2}{4w}s^2
\end{equation}
for every $\zeta\perp\nu$ and $s\in\mathbb R$.
\end{lemma}
\begin{proof}
The tangential, mixed and normal blocks of $\mathscr K_c[f]$ are
\[
 A+c\frac{p\otimes p}{w},\qquad
 b+c\frac{qp}{w},\qquad
 r+c\frac{q^2}{w}.
\]
The tangential block is at least $aI$, and the mixed block has
norm at most
\[
 J=L+c_+q_+R,
\]
while the normal block is bounded below by
$-L+c_-q_-^2/w$.  Therefore
\begin{align*}
 \mathscr K_c[f][\zeta+s\nu,\zeta+s\nu]
 &\ge a|\zeta|^2-2J|\zeta||s|
   +\left(\frac{c_-q_-^2}{w}-L\right)s^2\\
 &\ge \frac a2|\zeta|^2+
 \left(\frac{c_-q_-^2}{w}-L-\frac{2J^2}{a}\right)s^2.
\end{align*}
Choose $w_0$ so that
\[
 w_0\left(L+\frac{2J^2}{a}\right)
 \le\frac12c_-q_-^2.
\]
With this choice of $w_0$ (decreased if necessary), the last
coefficient is at least $c_-q_-^2/(4w)$, proving
\eqref{eq:K-lower}.
\end{proof}

\subsection{Convexity near the boundary}

The next proposition gives the general boundary argument used in
Corollary~\ref{cor:boundary-scale}.  It applies to both the square-root
and logarithmic transforms, with $f_D=U,\varphi_D$ and
$f_\alpha=u_\alpha,\varphi_\alpha$.
\begin{proposition}[From Dirichlet convexity to Robin convexity]\label{prop:transfer}
Let $\Omega$ be a bounded, uniformly convex, and $C^2$ domain.  Let
$f_D\in C^2(\overline\Omega)$ be negative in $\Omega$, vanish on
$\partial\Omega$, and satisfy $q_D:=(f_D)_\nu>0$ there.  Suppose
$f_\alpha\in C^2(\overline\Omega)$ is negative in $\Omega$ and
\[
 (f_\alpha)_\nu+\alpha f_\alpha=0\quad\text{on }\partial\Omega,
 \qquad \|f_\alpha-f_D\|_{C^2(\overline\Omega)}\longrightarrow0
          \quad\text{as }\alpha\to\infty.
\]
There exist $\alpha_0,\delta,c>0$, independent of $\alpha$, such that
$f_\alpha<0$ on $\overline\Omega$ and both estimates below hold
for $\alpha\ge\alpha_0$, $x=y-d\nu(y)$ with $0\le d<\delta$,
$\zeta\in T_y\partial\Omega$, and $s\in\mathbb R$:
\begin{align}
 D^2[-\sqrt{-f_\alpha}]
   [\zeta+s\nu,\zeta+s\nu]
 &\ge c\Bigl((d+\alpha^{-1})^{-1/2}|\zeta|^2
 +(d+\alpha^{-1})^{-3/2}s^2\Bigr),\label{eq:anisotropic-torsion}\\
 D^2[-\log(-f_\alpha)]
   [\zeta+s\nu,\zeta+s\nu]
 &\ge c\Bigl((d+\alpha^{-1})^{-1}|\zeta|^2
 +(d+\alpha^{-1})^{-2}s^2\Bigr).\label{eq:anisotropic-eigen}
\end{align}
These near-boundary estimates require no convexity assumption on a
Dirichlet transform.  If, in addition, a chosen transform
$\mathcal T(f)=-\sqrt{-f}$ or $\mathcal T(f)=-\log(-f)$ satisfies
$D^2\mathcal T(f_D)>0$ throughout $\Omega$, then there are
$\alpha_1\ge\alpha_0$ and $c_0>0$ such that
\[
 D^2\mathcal T(f_\alpha)\ge c_0I
       \quad\text{in }\Omega,\qquad \alpha\ge\alpha_1.
\]
The global assertion is made only for the transform satisfying this
additional hypothesis.
\end{proposition}

\begin{proof}
The first four steps prove the two boundary estimates using only
the geometric hypotheses.  In the fifth step, the additional
Dirichlet convexity assumption gives the corresponding global bound.

\smallskip
\noindent\emph{Step 1: tangential positivity and the normal slope.}
Uniform convexity and \eqref{eq:boundary-hessian} give
\[
 D^2f_D|_{T_y\partial\Omega}=q_D(y)\II_y
 \ge q_{D,*}\kappa_0I,
 \qquad q_{D,*}:=\min_{\partial\Omega}q_D>0.
\]
Fix a sufficiently small neighborhood of $\partial\Omega$.
By continuity and $C^2$ convergence, the tangential block $A_\alpha$
of $D^2f_\alpha$ then satisfies, for large $\alpha$,
\begin{equation}\label{eq:collar-tangent}
 A_\alpha\ge aI
\end{equation}
there, for a fixed $a>0$.  The $C^2$ convergence gives a uniform
bound for $D^2f_\alpha$.  Enlarging $L$ if necessary, the mixed and
normal blocks therefore satisfy $|b_\alpha|+|r_\alpha|\le L$.
The normal gradient component converges uniformly to that of $f_D$.
After shrinking the neighborhood if necessary,
\begin{equation}\label{eq:collar-normal}
 0<q_-\le Df_\alpha(x)\cdot\nu(y)\le q_+.
\end{equation}

\smallskip
\noindent\emph{Step 2: the boundary size of $-f_\alpha$.}
On $\partial\Omega$, the Robin condition yields
\[
 -f_\alpha(y)=\frac{(f_\alpha)_\nu(y)}{\alpha}.
\]
Along the normal segment $x=y-d\nu(y)$,
\[
 -f_\alpha(x)
 =-f_\alpha(y)+\int_0^d
   Df_\alpha(y-t\nu(y))\cdot\nu(y)\,dt.
\]
The bounds \eqref{eq:collar-normal} now give
\begin{equation}\label{eq:effective-distance}
 q_-(d+\alpha^{-1})
 \le -f_\alpha(x)
 \le q_+(d+\alpha^{-1})
\end{equation}
near the boundary.  Hence $-f_\alpha$ is comparable to $d+\alpha^{-1}$ near the
boundary.  This comparison will be used in the transformed-Hessian
bounds.

\smallskip
\noindent\emph{Step 3: the tangential gradient.}
Differentiating the Robin condition along a tangential vector
$\tau$ gives
\begin{equation}\label{eq:tangential-robin}
 \alpha\,\partial_\tau f_\alpha
 =-D^2f_\alpha[\tau,\nu]
   -Df_\alpha\cdot D_\tau\nu .
\end{equation}
Let $P_{T_y\partial\Omega}$ denote the orthogonal projection of
$\mathbb R^n$ onto the tangent space $T_y\partial\Omega$.  The uniform
Hessian and curvature bounds show that
\[
 |P_{T_y\partial\Omega}Df_\alpha(y)|\le C\alpha^{-1}
 \qquad(y\in\partial\Omega).
\]
For $x=y-d\nu(y)$, the Hessian bound gives
\[
 |P_{T_y\partial\Omega}Df_\alpha(x)|
 \le C(d+\alpha^{-1}).
\]
Together with \eqref{eq:effective-distance}, this gives
\begin{equation}\label{eq:tangential-gradient}
 |P_{T_y\partial\Omega}Df_\alpha(x)|
 \le R(-f_\alpha(x))
\end{equation}
for a fixed $R$.  In particular, the ratio $|p|/w$ in
Lemma~\ref{lem:schur} is uniformly bounded, so the mixed block of
$\mathscr K_c[f_\alpha]$ remains bounded even as $w=-f_\alpha$
becomes small.

\smallskip
\noindent\emph{Step 4: Schur positivity near the boundary.}
The preceding bounds meet the hypotheses of Lemma~\ref{lem:schur}
with $c_-=1/2$ and $c_+=1$, covering both $c=1/2$ and $c=1$.
Choose $\delta>0$ below the tubular radius and the inradius of
$\Omega$, with $q_+\delta\le w_0/4$, and then choose
$\alpha$ large enough that $q_+/\alpha\le w_0/4$.  For $d\le\delta$, \eqref{eq:effective-distance}
then implies
$w=-f_\alpha\le q_+(\delta+\alpha^{-1})\le w_0/2$.  Lemma~\ref{lem:schur}
therefore yields
\[
 \mathscr K_c[f_\alpha][\zeta+s\nu,\zeta+s\nu]
 \ge \frac a2|\zeta|^2+\frac{c_-q_-^2}{4(-f_\alpha)}s^2.
\]
With $w=-f_\alpha$, the two identities in
\eqref{eq:transform-formulas} give, respectively,
\begin{align*}
 D^2[-\sqrt{-f_\alpha}][\zeta+s\nu]
 &\ge \frac{a}{4w^{1/2}}|\zeta|^2
       +\frac{q_-^2}{16w^{3/2}}s^2,\\
 D^2[-\log(-f_\alpha)][\zeta+s\nu]
 &\ge \frac{a}{2w}|\zeta|^2
       +\frac{q_-^2}{8w^2}s^2.
\end{align*}
Since $q_-(d+\alpha^{-1})\le w\le q_+(d+\alpha^{-1})$
by \eqref{eq:effective-distance}, these are exactly
\eqref{eq:anisotropic-torsion}--\eqref{eq:anisotropic-eigen},
with constants independent of sufficiently large $\alpha$.

\smallskip
\noindent\emph{Step 5: the compact interior.}
Fix this value of $\delta>0$ and let
\[
 K_\delta=\{x\in\overline\Omega:d(x)\ge\delta\}.
\]
If $m_{f,\delta}=\min_{K_\delta}(-f_D)>0$, then the $C^0$
convergence gives $-f_\alpha\ge m_{f,\delta}/2$ on $K_\delta$
for sufficiently large $\alpha$.  In particular, all reciprocal
factors in the transformed Hessians are uniformly bounded there.
Let $\mathcal T(f)$ denote the relevant transform, namely
$-\sqrt{-f}$ or $-\log(-f)$.  By hypothesis,
$D^2\mathcal T(f_D)$ is positive definite on $K_\delta$.
Writing $\lambda_{\min}(A)$ for the smallest eigenvalue of a symmetric
matrix $A$, set
\[
 m_\delta
 =\min_{K_\delta}
   \lambda_{\min}\bigl(D^2\mathcal T(f_D)\bigr)>0.
\]
The map
\[
 (f,Df,D^2f)\longmapsto D^2\mathcal T(f)
\]
is smooth on the compact range of data over $K_\delta$.  Therefore
\[
 \|D^2\mathcal T(f_\alpha)
       -D^2\mathcal T(f_D)\|_{C^0(K_\delta)}
 \le C_\delta\|f_\alpha-f_D\|_{C^2(\overline\Omega)}.
\]
For sufficiently large $\alpha$, the right-hand side is at most
$m_\delta/2$.  The transformed Robin Hessian is then positive on
$K_\delta$.  The near-boundary lower bound and $m_\delta/2$ give a positive
minimum $c_0$ independent of $\alpha$.  Increasing $\alpha_1$ to
satisfy the preceding requirements proves the global assertion.
\end{proof}

\subsection{Completion of the main theorems}

For torsion, Section~\ref{sec:branches} establishes
\[
 \|u_\alpha-U\|_{C^{2,\gamma}}\le C\alpha^{-1}.
\]
The Dirichlet boundary condition and
Proposition~\ref{prop:dirichlet-input} give $U=0$ and $U_\nu>0$
on $\partial\Omega$, together with $D^2[-\sqrt{-U}]>0$ in
$\Omega$.  Proposition~\ref{prop:transfer} therefore yields
\[
 D^2[-\sqrt{-u_\alpha}]\ge c_TI.
\]
Near the boundary, \eqref{eq:anisotropic-torsion} gives the stronger
anisotropic bound.

Under \eqref{eq:Eg}, Section~\ref{sec:branches} establishes
\[
 \|\varphi_\alpha-\varphi_D\|_{C^{2,\gamma}}
   +|\lambda_\alpha-\lambda_D|
 \le C\alpha^{-1}.
\]
The Dirichlet boundary condition and
Proposition~\ref{prop:dirichlet-input} give $\varphi_D=0$ and
$(\varphi_D)_\nu>0$ on $\partial\Omega$, together with
$D^2[-\log(-\varphi_D)]>0$ in $\Omega$.  Applying the logarithmic case of
Proposition~\ref{prop:transfer} yields
\[
 D^2[-\log(-\varphi_\alpha)]\ge c_EI.
\]
The estimate \eqref{eq:anisotropic-eigen} also holds near the boundary.  Multiplying an
eigenfunction by a positive constant adds only a constant to
$-\log(-\varphi_\alpha)$, so this Hessian conclusion is independent of
the normalization.  The final thresholds $\alpha_T$ and $\alpha_E$ are taken to be
the maxima of the respective lower bounds required by the linear
estimates, the contraction arguments, the sign and admissibility
checks, and the boundary-transfer step.  Combined with the existence,
uniqueness, admissibility, and $O(\alpha^{-1})$ convergence established
in Section~\ref{sec:branches}, the bounds above prove
Theorems~\ref{thm:torsion-main} and~\ref{thm:eigen-main}.
Finally, \eqref{eq:effective-distance} and
\eqref{eq:anisotropic-torsion}--\eqref{eq:anisotropic-eigen}
prove Corollary~\ref{cor:boundary-scale}.

\section{Consequences and remarks}\label{sec:remarks}

We record a geometric consequence and the optimality of the convergence
rate, verify the spectral boundary hypothesis on balls, and discuss
the dependence of the constants.

\begin{corollary}[Convex superlevel sets]\label{cor:superlevel}
For $\alpha$ above the respective thresholds in the main theorems,
every nonempty superlevel set of $-u_\alpha$ and, under \eqref{eq:Eg},
of $-\varphi_\alpha$ is convex.  Indeed, $\sqrt{-u_\alpha}$ and $\log(-\varphi_\alpha)$
are strictly concave, and both the square root and logarithm are
increasing on the positive half-line.
\end{corollary}

\begin{remark}[Optimal scale of the convergence rate]\label{rem:sharp-rate}
The convergence rate $O(\alpha^{-1})$ for the Robin solution
functions is sharp already in $C^0$.  The Robin boundary
identities and the $C^1$ convergence imply
\[
 \alpha u_\alpha|_{\partial\Omega}\longrightarrow-U_\nu,
 \qquad
 \alpha\varphi_\alpha|_{\partial\Omega}
        \longrightarrow-(\varphi_D)_\nu
 \quad\text{in }C^0(\partial\Omega).
\]
Since $U_\nu\ge q_{T,*}>0$ and
$(\varphi_D)_\nu\ge q_{E,*}>0$, it follows that
\[
 \|u_\alpha-U\|_{C^0}\ge c\alpha^{-1},\qquad
 \|\varphi_\alpha-\varphi_D\|_{C^0}\ge c\alpha^{-1}
\]
for sufficiently large $\alpha$ (the second estimate under
\eqref{eq:Eg}).  Thus the scale in
Corollary~\ref{cor:boundary-scale} is both the Dirichlet
convergence rate of the solutions and their natural boundary
length scale.
\end{remark}

\begin{proposition}[The spectral boundary hypothesis on balls]
\label{prop:Eg-ball}
Let $B_R\subset\mathbb R^n$, $n\ge3$, be a ball.  The normalized
negative admissible Dirichlet $2$-Hessian eigenfunction on $B_R$
belongs to $C^\infty(\overline{B_R})$.  Consequently, \eqref{eq:Eg} holds on every ball for every
$0<\gamma<1$.
\end{proposition}
\begin{proof}
Let $(\lambda_D,\varphi_D)$ be the normalized Dirichlet eigenpair on
$B_R$.  Rotational invariance and uniqueness up to scale imply that
$\varphi_D(x)=\psi(r)$, $r=|x|$.  Wang's eigenfunction theory \cite[Theorem~6.2]{Wang09} yields
$\varphi_D\in C^{1,1}(\overline{B_R})$.  Interior elliptic
regularity then gives $\varphi_D\in C^\infty(B_R)$, while the Hopf
boundary lemma gives
\[
 \psi(R)=0,\qquad \psi'(R)>0.
\]
The eigenvalues of $D^2\varphi_D$ are $\psi''$ in the radial direction
and $\psi'/r$ with multiplicity $n-1$.  Consequently, for almost every $r$ near $R$,
\begin{equation}\label{eq:radial-eigen-ode}
 (n-1)\frac{\psi'}r\psi''
 +\binom{n-1}{2}\left(\frac{\psi'}r\right)^2
 =\lambda_D\psi^2.
\end{equation}
Since $\psi'$ is continuous and positive at $R$, it is bounded away from
zero on $[R-\varepsilon,R]$ for some $\varepsilon>0$.  Solving
\eqref{eq:radial-eigen-ode} for the second derivative gives
\begin{equation}\label{eq:radial-eigen-second}
 \psi''
 =\frac{r\lambda_D\psi^2}{(n-1)\psi'}
 -\frac{n-2}{2r}\psi'.
\end{equation}
The right-hand side is continuous on $[R-\varepsilon,R]$.
Because $\psi\in C^{1,1}$, its derivative $\psi'$ is absolutely
continuous.  Equation~\eqref{eq:radial-eigen-second} identifies the
almost-everywhere derivative of $\psi'$ with this continuous function.  Hence $\psi\in C^2$ up to $r=R$.
Because $\psi'(R)>0$, the right-hand side of
\eqref{eq:radial-eigen-second} is a smooth function of
$(r,\psi,\psi')$ near $r=R$.  Standard regularity for ordinary differential equations now yields
$\psi\in C^\infty([R-\varepsilon,R])$.  Together with the interior
regularity, this proves
$\varphi_D\in C^\infty(\overline{B_R})$.  In particular,
\[
 \psi''(R)=-\frac{n-2}{2R}\psi'(R).
\]
\end{proof}

\begin{remark}[Domains satisfying \texorpdfstring{$(E_\gamma)$}{E-gamma}]
\label{rem:Eg-domains}
By Proposition~\ref{prop:Eg-ball}, balls satisfy \eqref{eq:Eg}
automatically.  The domain perturbations in
Lemma~\ref{lem:spectral-approx} approximate a fixed Dirichlet pair
already satisfying \eqref{eq:Eg}; they do not establish that
hypothesis on an arbitrary uniformly convex $C^{3,\gamma}$ domain.
The latter remains a separate boundary-regularity question.
\end{remark}

\begin{remark}[Dependence and uniformity of the constants]
\label{rem:constant-dependence}
For a fixed domain and fixed normalized Dirichlet data, all constants
and parameter thresholds above are independent of $\alpha$.
For varying domains, separate the analytic estimates from the global
convexity conclusion.  Assume a common finite boundary atlas with
uniform chart radii and $C^{3,\gamma}$ bounds, a diameter bound, a
positive lower curvature bound, and common bounds for
\[
 \|U\|_{C^{3,\gamma}},\quad q_{T,*}^{-1},\quad
 \left(\min_{\overline\Omega}
         \lambda_{\min}(T_1(D^2U))\right)^{-1}.
\]
In the spectral case assume corresponding bounds for
$\|\varphi_D\|_{C^{3,\gamma}}$, $q_{E,*}^{-1}$, the ellipticity
constants of $a_E$, upper and positive lower bounds for $\lambda_D$,
and the inverse norm of \eqref{eq:dir-aug-map} at $\eta=\gamma$.
The linear estimates, contraction argument, and the near-boundary
bounds are uniform under these controls.  A fixed boundary width
$\delta$ may be chosen from these data.

For $\delta>0$, write
$K_\delta(\Omega)=\{x\in\overline\Omega:d(x)\ge\delta\}$.
When this set is nonempty, the global Hessian lower bounds
additionally use the compact-interior margins
\begin{equation}\label{eq:uniform-convexity-margins}
 \begin{split}
 m_{T,\delta}&=\min_{x\in K_\delta(\Omega)}
       \lambda_{\min}\bigl(D^2[-\sqrt{-U}](x)\bigr),\\
 m_{E,\delta}&=\min_{x\in K_\delta(\Omega)}
       \lambda_{\min}\bigl(D^2[-\log(-\varphi_D)](x)\bigr).
 \end{split}
\end{equation}
These numbers are positive for each fixed domain by
Proposition~\ref{prop:dirichlet-input}.  A common positive lower
bound for the relevant margins makes the global constants
$c_T,c_E$ and the corresponding concavity thresholds uniform over
a family as well.  Indeed, in Step~5 of
Proposition~\ref{prop:transfer} the perturbation error is
$C_\delta\alpha^{-1}$ and must be less than half of the relevant
margin.  The reciprocal factors in these transformed Hessians are also
uniformly controlled on $K_\delta(\Omega)$.  Indeed,
$2\sigma_2(D^2U)=(\Delta U)^2-|D^2U|^2=2$ and $\Delta U>0$
give $\Delta U\ge\sqrt2$, while
\eqref{eq:reference-trace-gap} gives
$\Delta\varphi_D\ge n\vartheta_E/(n-1)$.
More generally, if $f\le0$ and $\Delta f\ge c_*>0$,
then comparison on $B_\delta(x)\subset\Omega$, for
$x\in K_\delta(\Omega)$, with the function
$c_*(|y-x|^2-\delta^2)/(2n)$ gives
\[
 -f(x)\ge\frac{c_*\delta^2}{2n}.
\]
Consequently,
\[
 -U(x)\ge\frac{\sqrt2}{2n}\delta^2,\qquad
 -\varphi_D(x)\ge\frac{\vartheta_E}{2(n-1)}\delta^2
 \quad(x\in K_\delta(\Omega)).
\]
This controls the reciprocal factors by the stated geometric and
ellipticity bounds.  If $K_\delta(\Omega)$ is empty, the boundary
estimates already cover the whole domain.  We do not infer a uniform
interior convexity margin solely from strict convexity on each domain.
\end{remark}

\section*{Acknowledgments}
The authors would like to thank Professor Xi-Nan Ma for his encouragement and suggestions on this subject. This work was supported by the Program for Young Talents of Basic Research in University of Heilongjiang Province (Grant No. YQJH2025116). The authors acknowledge the use of AI tools. All mathematical statements and proofs were independently verified by the authors, who take full responsibility for the content of the manuscript.

\end{document}